\documentclass[12pt,a4paper]{article}
\usepackage[top=3cm, bottom=2cm, left=3.cm, right=3.cm]{geometry}
\usepackage{lastpage}
\usepackage{rotating}
\usepackage{comment}
\usepackage{pdflscape}
\usepackage{booktabs}
\usepackage{longtable}
\usepackage{adjustbox}
\usepackage{standalone}
\usepackage{amsmath, amsthm}
\usepackage[dvipsnames]{xcolor}
\usepackage[framemethod=TikZ]{mdframed}
\usepackage{enumerate}
\usepackage[shortlabels]{enumitem}
\usepackage{fancyhdr}
\usepackage{indentfirst}
\usepackage{authblk}

\usepackage{listings}
\usepackage{sectsty}
\usepackage{thmtools}
\usepackage{subcaption}
\usepackage[dvipsnames]{xcolor}
\usepackage{tabularx}
\usepackage{shadethm}
\usepackage{hyperref}
\usepackage{setspace}
\usepackage{algorithmic}
\usepackage{tikz}
\usetikzlibrary{positioning}
\usepackage{graphicx}
\usepackage{mathtools,amssymb}
\usepackage{booktabs}
\usepackage{multirow}
\usepackage{adjustbox}
\usepackage{siunitx}
\usepackage{setspace}
\usepackage[ruled,vlined]{algorithm2e}
\SetAlFnt{\small}
\SetAlCapFnt{\small}
\SetAlCapNameFnt{\small}
\SetKwInput{KwIn}{Input}
\SetKwInput{KwOut}{Output}
\DontPrintSemicolon
\usepackage{xcolor}

\SetCommentSty{mycommfont}
\RestyleAlgo{ruled} 
\usepackage{natbib}
\setcitestyle{authoryear,open={(},close={)}}
\usepackage{titlesec}
\usepackage{textcomp}
\titlelabel{\thetitle.\quad}
\usepackage[dvipsnames]{xcolor}
\usepackage{siunitx}  % Include siunitx package for rounding and formatting numbers
\usepackage{amsmath}

\usepackage{fancyvrb}
\usepackage{natbib}
\makeatletter
\newcommand\mathcircled[1]{%
	\mathpalette\@mathcircled{#1}%
}
\newcommand\@mathcircled[2]{%
	\tikz[baseline=(math.base)] \node[draw,circle,inner sep=1pt] (math) {$\m@th#1#2$};%
}
\usepackage{fancyhdr}
\usepackage{blindtext}
\makeatother
\usepackage{amsmath} % for align environment
\hypersetup{
	colorlinks=true,
	citecolor=red,        % \citep, \citet, \citealp
	linkcolor=blue,       % \ref, \eqref, page and ToC links
	urlcolor=blue,        % \url, \href
	filecolor=magenta,    % links to local files
	bookmarksnumbered=true,
	pdftitle={A scenario-cluster-based and enhanced progressive hedging algorithm for the two-stage stochastic quadratic knapsack problem},
	pdfauthor={Ibrahim Dan Dije, Franklin Djeumou Fomeni, Leandro C. Coelho, Janosch Ortmann},
}
\usepackage{thmtools}
\usepackage[framemethod=TikZ]{mdframed}
\usepackage{tikz}

\mdfsetup{skipabove=\topskip,skipbelow=\topskip}
\newrobustcmd\ExampleText{%
	An \textit{inhomogeneous linear} differential equation has the form
	\begin{align}
		L[v ] = f,
	\end{align}
	where $L$ is a linear differential operator, $v$ is the dependent
	variable and $f$ is a given non−zero function of the independent
	variables alone.
}
\mdfdefinestyle{theoremstyle}{%
	linecolor=black,linewidth=1pt,%
	frametitlerule=true,%
	frametitlebackgroundcolor=gray!20,
	innertopmargin=\topskip,
}
\mdtheorem[style=theoremstyle]{Problem}{Problem}
\definecolor{codegreen}{rgb}{0,0.6,0}
\definecolor{codegray}{rgb}{0.5,0.5,0.5}
\definecolor{codepurple}{rgb}{0.58,0,0.82}
\definecolor{backcolour}{rgb}{0.95,0.95,0.92}

\lstdefinestyle{mystyle}{
	backgroundcolor=\color{backcolour},   
	commentstyle=\color{codegreen},
	keywordstyle=\color{magenta},
	numberstyle=\tiny\color{codegray},
	stringstyle=\color{codepurple},
	basicstyle=\ttfamily\footnotesize,
	breakatwhitespace=false,         
	breaklines=true,                 
	captionpos=b,                    
	keepspaces=true,                 
	numbers=left,                    
	numbersep=5pt,                  
	showspaces=false,                
	showstringspaces=false,
	showtabs=false,                  
	tabsize=2
}
\usetikzlibrary{positioning}

\usepackage{float}      % for [H]
\usepackage{needspace}  % for \Needspace
\catcode`\<=\active \def<{
	\fontencoding{T1}\selectfont\symbol{60}\fontencoding{\encodingdefault}}
\catcode`\>=\active \def>{
	\fontencoding{T1}\selectfont\symbol{62}\fontencoding{\encodingdefault}}
\catcode`\<=\active \def<{
	\fontencoding{T1}\selectfont\symbol{60}\fontencoding{\encodingdefault}}

\usepackage{amssymb}
\begin{document}

\title{A scenario-cluster-based and enhanced progressive hedging algorithm for the two-stage stochastic quadratic knapsack problem}

\author[1]{Ibrahim Dan Dije}
\author[1]{Franklin Djeumou Fomeni}
\author[2]{Leandro C. Coelho}
\author[3]{Janosch Ortmann}

\affil[1]{GERAD, CIRRELT \& Department of Analytics, Operations and Information Technology, Université du Québec à Montréal, 1250 Rue Sanguinet, Montréal, Québec H2X 3E7, Canada}

\affil[2]{GERAD, CIRRELT \& Department of Operations and Decision Systems, Université Laval, Québec, Québec G1V 0A6, Canada}

\affil[3]{GERAD, CRM \& Department of Analytics, Operations and Information Technology, Université du Québec à Montréal, 1250 Rue Sanguinet, Montréal, Québec H2X 3E7, Canada}

\affil[ ]{\texttt{dan\_dije.ibrahim@courrier.uqam.ca}, \texttt{djeumou\_fomeni.franklin@uqam.ca}, \texttt{leandro.coelho@fsa.ulaval.ca}, \texttt{ortmann.janosch@uqam.ca}}

\date{}

\maketitle

\begin{abstract}
\small
This paper introduces a two-stage stochastic quadratic knapsack problem (TSSQKP) with uncertainty in both profits and weights. To overcome the computational difficulties arising from the combined stochastic, binary, and quadratic structure, we propose a solution framework that integrates an adaptation of the progressive hedging algorithm with clustering-based scenario reduction for computing lower and upper bounds. The enhanced progressive hedging algorithm (EPHA) incorporates a rounding procedure and a dynamic penalty update strategy driven by the detection of oscillation and stagnation patterns. In parallel, we propose a clustering framework based on opportunity-cost distances between scenarios, from which we derive the medoid lower bound ($\mathrm{LB}_{\mathrm{M}}$), the cluster lower bound ($\mathrm{LB}_{\mathrm{C}}$), and the cluster upper bound ($\mathrm{UB}_{\mathrm{C}}$). This decomposition breaks the original problem into smaller subproblems, thereby reducing computational time and memory requirements while approximating the solution of the full-scenario problem. To further improve computational efficiency, EPHA is also used to solve the cluster subproblems heuristically; in this case, the resulting upper bound is referred to as the estimated cluster upper bound ($\mathrm{EUB}_{\mathrm{C}}$), and the corresponding gaps are interpreted as estimated optimality gaps. Computational experiments on \(800\) generated TSSQKP instances show that the proposed framework yields smaller gaps than CPLEX under comparable computational conditions. In the EPHA-based clustering framework, the average estimated gaps associated with $\mathrm{LB}_{\mathrm{C}}$, $\mathrm{LB}_{\mathrm{EPHA}}$, and $\mathrm{LB}_{\mathrm{M}}$ are \(1.65\%\), \(1.88\%\), and \(2.41\%\), respectively, compared with an average optimality gap of \(35.86\%\) for CPLEX. The computation of $\mathrm{LB}_{\mathrm{M}}$ is the fastest, with an average execution time of \(420.88\) seconds, while EPHA provides a favorable compromise between solution quality and computational effort, reducing the average computational time associated with the cluster bounds by approximately \(52.7\%\), from \(5156.79\) to \(2439.10\) seconds, while increasing the average estimated gap by only \(0.23\) percentage points.
\end{abstract}


\section{Introduction}
\label{Intro}
The Quadratic Knapsack Problem (QKP) is a well-known NP-hard combinatorial optimization problem that extends the classical knapsack problem by incorporating pairwise interactions between selected items through quadratic profit terms. Due to its ability to capture complementarities and interaction effects among decisions, the QKP has numerous applications in portfolio optimization, transportation planning, telecommunications, production systems, and resource allocation. The practical relevance and computational difficulty of the QKP have motivated extensive research on exact methods, relaxations, decomposition techniques, and heuristics \citep{pisinger2007quadratic,cacchiani2022knapsack,galli2025quadratic}.

%Despite the extensive literature devoted to deterministic QKP formulations, 
Many practical applications involve random parameters that cannot be accurately represented using deterministic models \citep{cheng2014distributionally}. In finance, transportation, logistics, and production planning, profits, weights, and capacities are subject to many sources of uncertainty such as changing market prices, fluctuating item performance, measurement errors, and inaccurate estimates of size, length, or weight. Ignoring this uncertainty may produce solutions that become inefficient or infeasible once the actual data are realized. Stochastic optimization provides a more realistic framework for sequential decision-making in the presence of incomplete information.

Two-stage stochastic optimization provides a natural framework for modeling sequential decisions under uncertainty. In this setting, first-stage decisions are determined before uncertainty realization, whereas second-stage recourse actions are performed after additional information becomes available. Such a structure is particularly suitable for knapsack-type applications in which initial selections may require subsequent adjustments once stochastic parameters are observed. Nevertheless, combining binary recourse variables with quadratic objective functions yields highly challenging stochastic binary quadratic optimization problems. Some well-motivated works related to nonlinear two-stage stochastic problem formulations can be found in \cite{mehrotra2009decomposition, li2018improved, li2019finite}.

Motivated by these observations, we introduce a two-stage stochastic quadratic knapsack problem (TSSQKP) with recourse decisions. The proposed model considers uncertainty in second-stage profits and weights while allowing corrective actions through item addition and removal decisions after information is revealed. The resulting formulation leads to a large-scale stochastic binary quadratic optimization problem whose size rapidly increases with the number of scenarios. To address these scalability issues, we investigate decomposition-based and clustering-based methodologies inspired by recent advances in stochastic optimization \citep{hewitt2022decision,keutchayan2023problem}. In particular, we combine the enhanced progressive hedging algorith (PHA) with decision-based scenario clustering techniques to derive efficient lower and upper bounds for the TSSQKP.

A further major challenge arises from the large number of scenarios required to accurately approximate probability distributions. Although larger scenario sets generally improve the quality of approximations, they also increase computational effort and may render stochastic formulations intractable in practice \citep{hewitt2022decision}. The central motivation of this study is to address the computational complexity arising from the nonlinear recourse structure, together with the computational challenges arising from the increasing number of scenarios.

The main contributions of this work are as follows. First, we introduce a two-stage stochastic extension of the QKP that incorporates binary recourse decisions under uncertain profits and weights, which is different from the formulations proposed in \cite{lisser2010stochastic,gaivoronski2011knapsack}. Second, we propose an enhancement to the well-known PHA \citep{rockafellar1991scenarios}, combined with decision-based scenario clustering techniques \citep{hewitt2022decision}, to derive efficient lower and upper bounding procedures for the proposed model. Finally, we demonstrate the effectiveness of the proposed framework through computational experiments on TSSQKP instances with varying numbers of items, scenario sizes, and profit matrix densities.

The remainder of this paper is organized as follows. Section \ref{Lit_Review} reviews the related literature on deterministic and stochastic QKP as well as scenario reduction and clustering techniques. Section \ref{Prop_Model} presents the proposed two-stage stochastic quadratic formulation and its linearized scenario-based counterpart. Sections \ref{Solution_Meth_PHA} and \ref{Solution_Meth_SRM} describe the proposed solution methodologies, including the enhanced progressive hedging algorithm framework and the scenario reduction bounding methods. Section \ref{Num_Ex} reports the computational experiments and discusses the numerical results. Finally, Section \ref{Concl} concludes, discusses the results of the paper and outlines future research directions.

\section{Literature review}
\label{Lit_Review}
The consideration of uncertainty is already well studied in the literature of the linear knapsack problem. Indeed, \cite{kosuch2011two} investigate a two-stage stochastic knapsack problem with independently normally distributed item weights. They propose upper and lower bounding procedures combined with a branch-and-bound algorithm to solve the problem under continuous uncertainty, without relying on a finite scenario discretization. They also discuss the possibility of applying a sample average approximation approach by discretizing the continuous sample space into a finite number of sampled scenarios, while noting that such an approach becomes impractical when the required sample size is large, even in the linear stochastic knapsack case. 
\cite{cohn1998stochastic} consider a stochastic knapsack problem with independent normally distributed item weights in the context of robust transportation planning and develop heuristic solution approaches, whereas \cite{dean2008approximating} investigate adaptive and non-adaptive policies and provide constant-factor approximation guarantees that quantify the benefit of adaptivity. However, neither study adopts a two-stage, scenario-based stochastic programming formulation.

On the other hand, the QKP, which generalizes the linear knapsack problem, has primarily been studied from a deterministic perspective. Indeed, several algorithms have been developed for the deterministic QKP, including constructive and dynamic programming heuristics \citep{julstrom2005greedy,fomeni2014dynamic,djeumou2023lifted,EliassDP,hochbaum2025fast}, metaheuristic approaches based on genetic algorithms, Greedy Randomized Adaptive Search Procedures (GRASP), tabu search, and iterated search \citep{julstrom2005greedy,yang2013effective,chen2017iterated}, as well as exact methods based on branch-and-bound, Lagrangian decomposition, aggressive reduction, and $t$-linearization techniques \citep{gallo1980quadratic,cap1999,billionnet2004exact,pisinger2007solution,rodrigues2012exact,fennich2025tight}.

The stochastic extension of the QKP has comparatively received litle attention. This is likely because combining the nonlinear nature of the QKP with recourse decisions significantly increases the complexity of the problem.

Among the pioneering contributions in this direction, \cite{lisser2010stochastic,gaivoronski2011knapsack} studied two-stage stochastic quadratic knapsack models with probability constraints and recourse mechanisms. To address the resulting stochastic quadratic binary programs, the authors developed semidefinite programming (SDP) relaxations and RLT-1 formulations to obtain tight upper bounds. Their study also highlighted the high computational cost induced by the number of scenarios and items size. Specifically, \cite{lisser2010stochastic}, limited their tests to $20$ items and $20$ scenarios.

Distributionally robust variants of the QKP have also been explored. In \cite{cheng2014distributionally}, a distributionally robust stochastic QKP with chance-constrained knapsack constraints is considered,
where only partial information on the uncertain data is available, such as moments, joint support, and possibly independence assumptions. The study focuses on semidefinite programming (SDP) relaxations of the binary formulation and develops tractable upper and lower bounding approaches for the resulting model.

A recurring difficulty in scenario-based stochastic optimization is the large number of scenarios required to accurately approximate probability distributions, which may render the resulting formulations intractable in practice \citep{hewitt2022decision}. To address this issue, several scenario reduction and clustering strategies have been proposed. Recently, scenario clustering techniques based on decision similarity have attracted increasing attention in stochastic optimization. These approaches group scenarios according to their impact on optimization decisions and solution structures. Such methodologies have been investigated in several contexts, including stochastic programming, scenario reduction, and decomposition-based optimization frameworks \citep{pflug2001scenario,dupacova2003scenario,hvattum2009using,bertsimas2009constructing,crainic2014scenario,hewitt2022decision}. In particular, opportunity-cost-based clustering methods have been shown to generate high-quality lower and upper bounds for stochastic optimization problems while reducing computational effort \citep{hewitt2022decision}. Among the methods developed for two-stage stochastic programs, the PHA has received particular attention, with its foundational theory, acceleration strategies, and successful applications well documented in the literature \citep{rockafellar1991scenarios,chun1995scenario,watson2011progressive,crainic2011progressive,gul2015progressive}.

Overall, the existing literature provides effective tools for either deterministic QKP or stochastic linear knapsack problems, but the combination of quadratic interactions, binary recourse decisions, and large scenario sets remains comparatively underexplored. This gap motivates the two-stage stochastic quadratic knapsack model and the clustering-based solution methodology developed in this paper.
\section{The two-stage stochastic quadratic knapsack problem}
\label{Prop_Model}
In this section, we present the formulations of the TSSQKP, including both the quadratic and the linearized forms. Our formulation differs from the stochastic quadratic knapsack models proposed in \cite{lisser2010stochastic,gaivoronski2011knapsack} in the treatment of knapsack capacity in the first and second stages. In \cite{lisser2010stochastic}, a probability constraint is imposed on the first-stage capacity, while a different random capacity is considered in the second stage. In contrast, \cite{gaivoronski2011knapsack} considers a deterministic capacity in the first stage and a potentially different random capacity in the second stage. In the latter formulation, the second-stage capacity constraint is imposed with a prescribed probability \(1-\alpha\), thereby allowing capacity violations for a subset of realizations whose total probability does not exceed \(\alpha\). In comparison, our formulation assumes that the physical capacity \(\beta\) remains unchanged between the two stages and imposes the second-stage capacity constraint for every scenario. Uncertainty is instead represented through the scenario-dependent item weights and profits rather than through changes in the available capacity.
\subsection{Stochastic quadratic formulation}
Let \(N=\{1,\ldots,n\}\) denote the set of items. The deterministic first-stage parameters are the symmetric quadratic profit matrix \(Q\), with nonnegative integer entries, the item-wight vector \(w\) with positive integer entries, and the integer knapsack capacity \(\beta>0\). Let \(\xi\) denote the random vector describing the uncertain second-stage parameters. For each realization of \(\xi\), \(T(\xi)\) and  \(R(\xi)\) denote the random profit and penalty matrices, respectively, while \(a_i(\xi)\) denotes the random weight of item \(i\in N\).

Let \(x\in\{0,1\}^n\) denote the first-stage decision variable, where \(x_i\) indicates the selection status of item \(i\in N\) in the first-stage. The TSSQKP is formulated as follows:
\begin{subequations}
\label{TSSQKP_Intro}
\begin{align}
\max_{x}\quad
& x^{T}Qx+
\mathbb{E}_{\xi}\!\left[\mathcal{Q}(x,\xi)\right]
\label{TSSQKP_Intro_Obj}
\\
\text{s.t.}\quad
& w^{T}x\leq \beta,
\label{TSSQKP_Intro_Capacity}
\\
& x_i\in\{0,1\},
\qquad \forall i\in N.
\label{TSSQKP_Intro_Binary}
\end{align}
\end{subequations}

For a given first-stage solution $x$ and a realization of the random vector $\xi$, the second-stage decision allows the initial selection to be adjusted once the uncertain item weights, profits and penalties are revealed. Two types of actions are possible in the second stage. First, an item selected in the first-stage can be removed from the knapsack. This decision is represented by the binary variable \(v_i\). Second, an item not selected in the first-stage can be added to the knapsack; this decision is represented by the binary variable \(u_i\). The parameter \(T(\xi)\) represents the uncertainty-dependent profit coefficients associated with the second-stage item selections, whereas \(R(\xi)\) represents the uncertainty-dependent penalty incurred when items selected in the first-stage are removed during the second stage. The second-stage recourse problem can then be formulated as follows:

\begin{subequations}
\label{TSSQKP_Recourse}
\begin{align}
\mathcal{Q}(x,\xi)=
\max_{u,v}\quad
& u^{T}T(\xi)u-v^{T}R(\xi)v
\label{TSSQKP_Recourse_Obj}
\\
\text{s.t.}\quad
& \sum_{i\in N}
a_i(\xi)\bigl(u_i+x_i-v_i\bigr)
\leq \beta,
\label{TSSQKP_Recourse_Capacity}
\\
& u_i\leq 1-x_i,
\qquad \forall i\in N,
\label{TSSQKP_Recourse_Addition}
\\
& v_i\leq x_i,
\qquad \forall i\in N,
\label{TSSQKP_Recourse_Removal}
\\
& u_i\in\{0,1\},
\qquad \forall i\in N,
\label{TSSQKP_Recourse_BinaryU}
\\
& v_i\in\{0,1\},
\qquad \forall i\in N.
\label{TSSQKP_Recourse_BinaryV}
\end{align}
\end{subequations}

The recourse function \(\mathcal{Q}(x,\xi)\) adjusts the initial first-stage solution after the uncertain profits and weights are revealed. The recourse process involves two corrective mechanisms. The binary variable $u_i$ indicates whether an item is added during the second stage, whereas $v_i$ indicates whether a previously selected item is removed after uncertainty realization. The second-stage capacity constraint guarantees feasibility of the updated solution, while the constraints \eqref{TSSQKP_Recourse_Addition} and \eqref{TSSQKP_Recourse_Removal} ensure consistency between first-stage and second-stage decisions. Specifically, an item may only be added if it was not selected initially, and it may only be removed if it was selected in the first-stage. 

This modelling choice is motivated by applications in which the available resource remains fixed, while information about item weights and profits is progressively revealed. It guarantees feasibility with respect to every scenario included in the model and avoids the need to specify a risk parameter \(\alpha\) \citep{shapiro2014lectures}. This requirement may produce a more conservative solution than a chance-constrained formulation, since no scenario-dependent capacity violation is permitted \citep{bertsimas2004price}. Nevertheless, it provides explicit scenario-wise recourse problems.

The proposed model combines two major sources of computational difficulty. First, the requirement that the first-stage binary decisions remain identical across all scenarios couples the scenario subproblems and prevents their independently obtained solutions from being directly combined. A change in a first-stage variable modifies the second stage decisions for every scenario, which makes the identification of a globally consistent first-stage solution challenging. Second, each scenario contains a binary quadratic recourse problem involving the variables \(u\) and \(v\). These quadratic terms introduce pairwise interactions between items and generally require either specialized mixed-integer quadratic optimization methods or a linearization with auxiliary variables and additional constraints. Consequently, the size of the deterministic equivalent grows rapidly with both the number of items and the number of scenarios.

\subsection{Linearized scenario-based formulation}
Using the finite scenario representation of the uncertainty as well as applying the McCormick linearization approach \citep{mccormick1976computability} on the quadratic terms appearing in \eqref{TSSQKP_Intro_Obj} and \eqref{TSSQKP_Recourse_Obj}, one obtains scenario-based deterministic equivalent model.
Let
\[
\mathcal{P}=\{(i,j)\in N\times N:i<j\}
\]
denote the set of unordered pairs of distinct items. The uncertainty associated with the second-stage profits and weights is represented by a finite scenario set \(\Omega\), where \(|\Omega|\) denotes the total number of scenarios. Each scenario \(\omega\in\Omega\) corresponds to a realization
\[
\xi^\omega=\left(a^\omega,T^\omega,R^\omega\right)
\]
of the random vector \(\xi\) and occurs with probability \(p_\omega>0\), where
\[
\sum_{\omega\in\Omega}p_\omega=1.
\]
When the scenarios are generated through sample-average approximation,
one may use \(p_\omega=1/|\Omega|\) for all \(\omega\in\Omega\). Then, the expected recourse value is approximated by
\[
\mathbb{E}_{\xi}\!\left[\mathcal{Q}(x,\xi)\right]
\approx
\sum_{\omega\in\Omega}p_\omega
\mathcal{Q}\left(x,\xi^\omega\right).
\]

For each scenario \(\omega\in\Omega\), \(a_i^\omega\) denotes the realized weight of item \(i\in N\), while \(T^\omega\) and \(R^\omega\) denote the scenario-dependent profit and removal cost matrices, respectively. The first-stage parameters are the profit matrix \(Q\), the deterministic weight vector \(w\), and the knapsack capacity \(\beta\).

The variables \(x_i\) indicate the items selected in the first-stage, whereas
\(u_i^\omega\) and \(v_i^\omega\) indicate, respectively, whether item \(i\)
is added to or removed from the knapsack under scenario \(\omega\). Before
presenting the deterministic equivalent formulation, we introduce the following auxiliary decision variables, which are used for linearization purpose.
\[
y_{ij}=x_ix_j,\qquad
z_{ij}^\omega=u_i^\omega u_j^\omega,\qquad
t_{ij}^\omega=v_i^\omega v_j^\omega,
\qquad (i,j)\in\mathcal{P},\ \omega\in\Omega.
\]
These binary products are linearized using the McCormick linearization technique \citep{mccormick1976computability}, then we have the linearized scenario-based formulation \eqref{TSSQKP_Linearise} which contains
\begin{equation}
 \frac{n(n+1)(2|\Omega|+1)}{2}
 \label{ComplexityVar}
\end{equation}

binary variables and
\begin{equation}
1+\frac{3n(n-1)}{2}
+|\Omega|\left(3n^2-n+1\right)
\label{ComplexityCon}
\end{equation}
linear constraints.
In the model \eqref{TSSQKP_Linearise}, for \(|\Omega|\) scenarios and $n$ items, both the number of variables and the number of constraints are of order \(\mathcal{O}(|\Omega|n^2)\) (see \eqref{ComplexityVar} and \eqref{ComplexityCon}).
\begin{center}
\begin{minipage}{0.95\textwidth}
\begingroup

% Compact vertical spacing
\setlength{\jot}{-3pt}
\setlength{\abovedisplayskip}{2pt}
\setlength{\belowdisplayskip}{2pt}
\setlength{\abovedisplayshortskip}{2pt}
\setlength{\belowdisplayshortskip}{2pt}

\begin{subequations}
\label{TSSQKP_Linearise}

%--------------------------------------------------------
% Objective function
%--------------------------------------------------------
\begin{equation}
\resizebox{0.98\linewidth}{!}{$
\displaystyle
\max_{x,y,u,v,z,t}\quad
\mathcal{F}(x;\Omega)
=
\sum_{i\in N} Q_{ii}x_i
+
2\sum_{(i,j)\in\mathcal{P}} Q_{ij}y_{ij}
+
\sum_{\omega\in\Omega}p_\omega
\left(
\sum_{i\in N}T_{ii}^{\omega}u_i^\omega
+
2\sum_{(i,j)\in\mathcal{P}}T_{ij}^{\omega}z_{ij}^\omega
-
\sum_{i\in N}R_{ii}^{\omega}v_i^\omega
-
2\sum_{(i,j)\in\mathcal{P}}R_{ij}^{\omega}t_{ij}^\omega
\right)
$}
\label{TSSQKP_Linearise_Obj}
\end{equation}

\vspace{-5pt}

%--------------------------------------------------------
% Constraints
%--------------------------------------------------------
\footnotesize

\begin{alignat}{2}
\text{s.t.}\quad
&
\sum_{i\in N} w_i x_i \leq \beta,
&&
\label{TSSQKP_Linearise_FS_Capacity}
\\
&
\sum_{i\in N}a_i^\omega
\left(
u_i^\omega+x_i-v_i^\omega
\right)
\leq \beta,
&\quad&
\forall \omega\in\Omega,
\label{TSSQKP_Linearise_SS_Capacity}
\\
&
y_{ij}\leq x_i,
&\quad&
\forall (i,j)\in\mathcal{P},
\label{TSSQKP_Linearise_Y1}
\\
&
y_{ij}\leq x_j,
&\quad&
\forall (i,j)\in\mathcal{P},
\label{TSSQKP_Linearise_Y2}
\\
&
y_{ij}\geq x_i+x_j-1,
&\quad&
\forall (i,j)\in\mathcal{P},
\label{TSSQKP_Linearise_Y3}
\\
&
u_i^\omega\leq 1-x_i,
&\quad&
\forall i\in N,\ \omega\in\Omega,
\label{TSSQKP_Linearise_Addition}
\\
&
v_i^\omega\leq x_i,
&\quad&
\forall i\in N,\ \omega\in\Omega,
\label{TSSQKP_Linearise_Removal}
\\
&
z_{ij}^\omega\leq u_i^\omega,
&\quad&
\forall (i,j)\in\mathcal{P},\ \omega\in\Omega,
\label{TSSQKP_Linearise_Z1}
\\
&
z_{ij}^\omega\leq u_j^\omega,
&\quad&
\forall (i,j)\in\mathcal{P},\ \omega\in\Omega,
\label{TSSQKP_Linearise_Z2}
\\
&
z_{ij}^\omega
\geq u_i^\omega+u_j^\omega-1,
&\quad&
\forall (i,j)\in\mathcal{P},\ \omega\in\Omega,
\label{TSSQKP_Linearise_Z3}
\\
&
t_{ij}^\omega\leq v_i^\omega,
&\quad&
\forall (i,j)\in\mathcal{P},\ \omega\in\Omega,
\label{TSSQKP_Linearise_T1}
\\
&
t_{ij}^\omega\leq v_j^\omega,
&\quad&
\forall (i,j)\in\mathcal{P},\ \omega\in\Omega,
\label{TSSQKP_Linearise_T2}
\\
&
t_{ij}^\omega
\geq v_i^\omega+v_j^\omega-1,
&\quad&
\forall (i,j)\in\mathcal{P},\ \omega\in\Omega,
\label{TSSQKP_Linearise_T3}
\\
&
x_i,u_i^\omega,v_i^\omega\in\{0,1\},
&\quad&
\forall i\in N,\ \omega\in\Omega,
\label{TSSQKP_Linearise_Binary1}
\\
&
y_{ij},z_{ij}^\omega,t_{ij}^\omega\in\{0,1\},
&\quad&
\forall (i,j)\in\mathcal{P},\ \omega\in\Omega.
\label{TSSQKP_Linearise_Binary2}
\end{alignat}

\end{subequations}

\endgroup
\end{minipage}
\end{center}

For the remainder of the paper, \(\mathcal{F}(x;\Omega)\) denotes the TSSQKP objective associated with the complete scenario set \(\Omega\), while
\(f(x;\omega)\) denotes its single-scenario counterpart for \(\omega\in\Omega\), with $\mathcal{F}(x;\Omega)=\frac{1}{|\Omega|}\sum_{\omega\in\Omega}f(x;\omega)$. Moreover, \(\mathcal{X}_{\Omega}\) and
\(\mathcal{X}_{\omega}\) denote the feasible regions associated with the
complete scenario-based formulation and the corresponding single-scenario
formulation, respectively.

\section{Progressive hedging method}
\label{Solution_Meth_PHA}
In this section, we propose an enhancement of the well-studied PHA \citep{rockafellar1991scenarios,crainic2011progressive,gade2016obtaining}, which we adapt using a rounding strategy and stagnation detection heuristics. Rounding is a well-established procedure in mixed-integer optimization that converts a fractional solution into an integer solution when specific conditions are satisfied \citep{achterberg2012rounding}. Stagnation-detection mechanisms are generally used to identify an absence of improvement and to adapt algorithmic parameters accordingly \citep{rajabi2023stagnation}. \cite{rajabi2023stagnation} use detection stagnation to automatically adjusts the mutation rate of evolutionary algorithms when they encounter local optima.

\subsection{Classical progressive hedging algorithm}
\label{ClassicalPHA}

The PHA is a powerful decomposition technique used to solve large-scale stochastic programming problems, particularly those with a block-separable structure. It was first introduced by \cite{rockafellar1991scenarios} and has since become a fundamental approach in stochastic optimization. PHA works by decomposing a large stochastic optimization problem into smaller, scenario-specific subproblems, which are then coordinated by means of an augmented Lagrangian function. The algorithm iteratively updates the Lagrangian multipliers and quadratic penalty terms which come directly from the augmented Lagrangian relaxation of the nonanticipativity constraints, ensuring that first-stage decisions are consistent across all scenarios. This iterative process converges to an optimal solution for convex problems and often provides good heuristic solutions for non-convex and mixed-integer problems. Its adaptability to various problem types, including those with binary variables, makes it a valuable tool in addressing the computational challenges of stochastic optimization.
For a deeper understanding of the PHA and its applications, especially in mixed-integer linear cases, the following works are highly relevant to the development and application of the PHA \citep{crainic2011progressive,gade2016obtaining,kaisermayer2021progressive}.

When applied to the TSSQKP, the PHA introduces a scenario-specific first-stage binary decision variable \(x_{\omega}\)
for each scenario \(\omega \in \Omega\). Each variable \(x_{\omega}\) is associated with the second-stage variables and constraints corresponding to scenario \(\omega\), which defines an independent single-scenario TSSQKP subproblem. Since the first-stage decisions must be made before the realized scenario is known, all scenario-specific first-stage solutions must satisfy the nonanticipativity conditions
\begin{equation}
x_{\omega}=x_{\omega'},
\qquad
\forall\,\omega,\omega' \in \Omega.
\label{PHA_nonanticipativity}
\end{equation}
Rather than imposing these equalities directly, at iteration $k$, PHA coordinates the scenario-related decisions through the probability-weighted consensus vector
\begin{equation}
\bar{x}^{(k)}
=
\sum_{\omega\in\Omega}
p_{\omega}x_{\omega}^{(k)},
\label{PHA_consensus}
\end{equation}
where \(x_{\omega}^{(k)}\) is the first-stage solution associated with scenario \(\omega\) at iteration \(k\).

At iteration \(k+1\), each scenario subproblem is solved independently by augmenting its original TSSQKP objective with a multiplier term and a quadratic penalty term:
\begin{equation}
\begin{aligned}
\max_{x
\in \mathcal{X}_{\omega}} \quad &
f
\left(
x;\omega
\right)
-
\left(\lambda_{\omega}^{(k)}\right)^{\top}
\left(
x-\bar{x}^{(k)}
\right)
-
\frac{\rho}{2}
\left\|
x-\bar{x}^{(k)}
\right\|_{2}^{2}
\end{aligned}
\label{PHA_penalized_subproblem}
\end{equation}
to obtain $x_{\omega}^{(k+1)}$, where \(\lambda_{\omega}^{(k)}\) is the multiplier vector associated with scenario \(\omega\), and \(\rho>0\) is the penalty parameter. The original objective favors decisions that perform well under scenario \(\omega\), whereas the augmented quadratic Lagrangian terms penalize deviations from the current consensus $\bar{x}^{(k)}$. After all scenario subproblems have been solved, the consensus vector and the multipliers are updated as follows:
\begin{align}
\bar{x}^{(k+1)}
&=
\sum_{\omega\in\Omega}
p_{\omega}x_{\omega}^{(k+1)},
\label{PHA_consensus_update}
\\
\lambda_{\omega}^{(k+1)}
&=
\lambda_{\omega}^{(k)}
+
\rho
\left(
x_{\omega}^{(k+1)}
-
\bar{x}^{(k+1)}
\right),
\qquad
\forall\,\omega\in\Omega.
\label{PHA_multiplier_update}
\end{align}
These updates progressively drive the scenario-specific first-stage decisions toward a common solution. The algorithm terminates when the maximum deviation between the scenario-specific solutions and the consensus is below a prescribed tolerance \(\epsilon>0\), that is,
\begin{equation}
\alpha^{(k+1)}
=
\sum_{\omega\in\Omega}
p_{\omega}
\left\|
x_{\omega}^{(k+1)}
-
\bar{x}^{(k+1)}
\right\|_{2}
\leq \epsilon.
\label{PHA_stopping_criterion}
\end{equation}
Algorithm \ref{PHA_Classic} provides a feasible solution for the TSSQKP. The quality of the solution depends on the choice of $\rho$, the maximum number of iterations $K$, and the convergence tolerance $\epsilon$ \citep{christiansen2023study}, which are standard PHA parameters.
\begin{algorithm}[H]

\SetKwInOut{KwIn}{Input}
\SetKwInOut{KwOut}{Output}

\scalebox{0.8}{%
\begin{minipage}{1.25\linewidth}

\KwIn{Penalty parameter $\rho>0$, maximum number of iterations $K$, tolerance $\epsilon>0$.}
\KwOut{A near-optimal first-stage solution $\bar{x}^{(k+1)}$.}

\tcc{--- Initialization step ---}
\ForEach{$\omega \in \Omega$ \tcp*[f]{in parallel}}{
	Solve:
	\[
	x_{\omega}^{*(0)}
	=
	\arg\max_{x\in \mathcal{X}_{\omega}}
	f(x;\omega)
	\]
}

$\bar{x}^{(0)}
\leftarrow
\sum_{\omega \in \Omega}
p_{\omega}\,x_{\omega}^{*(0)}$\;

$\lambda_{\omega}^{(0)}
\leftarrow
\rho
\big(
x_{\omega}^{*(0)}-\bar{x}^{(0)}
\big),
\quad
\forall \omega \in \Omega$\;

\tcc{--- Main PHA loop ---}
\For{$k=0,\ldots,K-1$}{

	\ForEach{$\omega \in \Omega$ \tcp*[f]{in parallel}}{

		Solve the augmented subproblem:
		\begin{equation}
		x_{\omega}^{*(k+1)}
		=
		\arg\max_{x\in \mathcal{X}_{\omega}}
		\Big(
		f(x;\omega)
		-
		(\lambda_{\omega}^{(k)})^\top x
		-
		\tfrac{\rho}{2}
		\sum_{i\in N}
		\big(1-2\bar{x}_i^{(k)}\big)x_i
		\Big)
		\label{AugmL_C}
		\end{equation}
	}

	$\bar{x}^{(k+1)}
	\leftarrow
	\sum_{\omega \in \Omega}
	p_{\omega}\,x_{\omega}^{*(k+1)}$\;

	$\lambda_{\omega}^{(k+1)}
	\leftarrow
	\lambda_{\omega}^{(k)}
	+
	\rho
	\big(
	x_{\omega}^{*(k+1)}
	-
	\bar{x}^{(k+1)}
	\big),
	\quad
	\forall \omega \in \Omega$\;

	$\alpha^{(k+1)}
	\leftarrow
	\sum_{\omega \in \Omega}
	p_{\omega}\,
	\|
	x_{\omega}^{*(k+1)}
	-
	\bar{x}^{(k+1)}
	\|_{2}$\;

	\If{$\alpha^{(k+1)} \le \epsilon$}{
		\textbf{break}\;
	}
}

\Return $\bar{x}^{(k+1)}$\;

\end{minipage}%
}

\caption{The PHA for the TSSQKP}
\label{PHA_Classic}

\end{algorithm}

The linearized equation \eqref{AugmL_C} in the algorithm at iteration $k$ is the augmented Lagrangian given by:

\begin{equation}
	L^k(x,\rho,\lambda;\omega) = f(x;\omega)
	- \lambda_\omega^{T}(x - \bar{x}^k)
	- \frac{\rho}{2}\|x - \bar{x}^k\|_2^2.
	\label{AugmL_Non}
\end{equation}
Here $\|\cdot\|_2$ denotes the sum-of-squares ($\ell_2$) norm, which corresponds to the classical Progressive Hedging formulation. Alternative norms, such as the $\ell_1$ and $\ell_\infty$ norms, may also be considered for mixed-integer problems, as discussed by \cite{kaisermayer2021progressive}; however, they are not investigated in this work. The advantage of the $\ell_2$ norm is that in the present binary setting it admits an exact linear formulation. More precisely, optimizing the quadratic term $\frac{\rho}{2}\|x - \bar{x}^k\|_2^2$ is equivalent to optimizing
\begin{equation}
	\tfrac{\rho}{2}\left(\sum_{i \in N}(1 - 2\bar{x}_i^{(k)})x_i+\sum_{i \in N}\bar{x}_i^{2(k)}\right)
\end{equation}
because $x_i^2 = x_i$ for binary variables, and the constant terms involving $\bar{x}^k$ can be ignored. Let $x_{PHA}$ be a solution obtained from Algorithm \ref{PHA_Classic} feasible in $\mathcal{X}_{\Omega}$ and $x^{*}$ the optimal solution of the TSSQKP. Then we have

\begin{equation}
\mathcal{F}(x_{PHA};\Omega)\le \mathcal{F}(x^{*};\Omega).
\end{equation}

To improve Algorithm \ref{PHA_Classic}, we propose in the next section some enhancements.

\subsection{Enhanced progressive hedging algorithm}

In this section, we develop a simple rounding heuristic together with a procedure for detecting stagnation and oscillation in the convergence-rate sequence of the PHA. The consensus solution is rounded at each iteration, and whenever stagnation or oscillation is detected, the penalty parameter \(\rho\) is adaptively updated. It has been widely recognized that the penalty parameter \(\rho\) significantly affects the convergence behavior of the PHA, and several update strategies have been proposed in the literature, including fixed penalty schemes and adaptive update mechanisms \citep{hvattum2009using,crainic2011progressive,watson2011progressive}. Unlike these approaches, which typically update \(\rho\) according to the iterations, our strategy adapts \(\rho\) based on the detection of stagnation or oscillation patterns over a convergence sequence of iterations.

\subsubsection{Rounding heuristic}
Rounding is a commonly used heuristic in mixed-integer optimization for
converting fractional solution components into integer values when specified conditions are satisfied \citep{achterberg2012rounding}. Algorithm \ref{Rounding} presents a rounding heuristic that is integrated into the PHA framework described in Algorithm \ref{PHA_Classic} to construct a near binary consensus solution at each iteration.

More specifically, at each iteration \(k\), the consensus value \(\bar{x}_i^k\) for each \(i \in N\) is rounded to \(1\) if \(\bar{x}_i^k \geq \overline{\kappa}\), indicating that at least \(100\overline{\kappa}\%\) of the scenarios assign the value \(1\) to the corresponding first-stage variable. Similarly, \(\bar{x}_i^k\) is rounded to \(0\) if \(\bar{x}_i^k \leq \underline{\kappa}\), meaning that at most \(100\underline{\kappa}\%\) of the scenarios assign the value \(1\) to that variable. Otherwise, the value remains unchanged. The resulting rounded solution is then checked for feasibility. If the rounded solution is feasible, it is accepted as the updated consensus solution; otherwise, the original consensus vector is retained.

\begin{algorithm}[H]
\caption{Rounding Procedure}
\label{Rounding}

\scalebox{0.8}{%
\begin{minipage}{1.25\linewidth}

\KwIn{Consensus vector $\bar{x}^k$, upper threshold $\overline{\kappa}=0.9$, lower threshold $\underline{\kappa}=0.2$}
\KwOut{Rounded binary vector $\bar{\mu}^k$}

Initialize $\bar{\mu}^k \leftarrow \bar{x}^k$\;

\For{$i \in N$}{
	\eIf{$\bar{\mu}_i^k \leq \underline{\kappa}$}{
		$\bar{\mu}_i^k \leftarrow 0$\;
	}{
		\If{$\bar{\mu}_i^k \geq \overline{\kappa}$}{
			$\bar{\mu}_i^k \leftarrow 1$\;
		}
	}
}
\eIf{$\bar{\mu}^k$ is feasible}{
	\Return $\bar{\mu}^k$\;
}{
	\Return $\bar{x}^k$\;
}
\end{minipage}%
}
\end{algorithm}

Algorithm \ref{Rounding} may improve the convergence of PHA by converting components with strong scenario agreement into binary values. This reduces the disagreement between the scenario-specific first-stage decisions and the consensus vector, thereby helping the algorithm reach a common solution more rapidly.
\subsubsection{Stagnation detection heuristic}

As in many iterative optimization algorithms, the convergence rate of the PHA may exhibit oscillatory behavior or stagnation over successive iterations. Such behavior is often associated with the choice of the penalty parameter \(\rho\). Since \(\rho\) controls the penalty imposed on the deviation between the scenario-specific first-stage decisions \(x_{\omega}^{(k)}\) and the consensus solution \(\bar{x}^{(k)}\), an inappropriate value of \(\rho\) may slow down convergence or generate persistent oscillations.

To monitor the convergence behavior, let
\begin{equation}
\alpha^{(k)}
=
\sum_{\omega\in\Omega}
p_{\omega}
\left\|
x_{\omega}^{(k)}
-
\bar{x}^{(k)}
\right\|_{2}
\label{PHA_convergence_measure}
\end{equation}
denote the probability-weighted deviation of the scenario-specific first-stage solutions from the consensus solution at iteration \(k\). The sequence
\[
\mathcal{S}^{(k)}
=
\left(
\alpha_0,\alpha_1,\ldots,\alpha_k
\right).
\]
describes the evolution of the convergence measure up to iteration \(k\). Consensus is reached when \(\alpha_k\leq\epsilon\), where \(\epsilon>0\) is the stopping tolerance defined in Algorithm \ref{PHA_Classic}. When \(\alpha_k>\epsilon\), the recent values of \(\mathcal{S}^{(k)}\) are divided into \(m\) consecutive subsequences, each containing \(p\) convergence values. Therefore, at least \(mp\) convergence values are required before the stagnation-detection procedure can be applied. For each \(\ell=0,\ldots,m-1\), let \(\mathcal{W}_\ell\) denote the \(\ell\)-th most recent subsequence. The following parameters are used:

\begin{itemize}
 \item \(m\) is the number of consecutive subsequences considered;
 \item \(p\) is the number of convergence values contained in each
 subsequence;
 \item \(R_\ell\) is the range of subsequence \(\mathcal{W}_\ell\), which measures
 its oscillation amplitude;
 \item \(M_\ell\) is the midpoint of subsequence \(\mathcal{W}_\ell\), which
 represents the central level around which the convergence values oscillate;
 \item \(\varepsilon_R>0\) is the tolerance used to compare the ranges of
 the subsequences;
 \item \(\varepsilon_M>0\) is the tolerance used to compare the midpoints of the subsequences.
\end{itemize}

Oscillatory stagnation is detected when the ranges and midpoints remain nearly unchanged across the \(m\) subsequences. Algorithm
\ref{stagnation_detection} presents the proposed procedure.

\begin{algorithm}[H]
\caption{Stagnation detection}
\label{stagnation_detection}

\scalebox{0.7}{%
\begin{minipage}{1.25\linewidth}

\KwIn{Current iteration $k$, convergence sequence
$\mathcal{S}^{(k)}=(\alpha_0,\alpha_1,\ldots,\alpha_k)$,
number of subsequences $m\geq 2$, subsequence size $p$,
tolerances $\varepsilon_R>0$ and $\varepsilon_M>0$}

\KwOut{Boolean value indicating whether oscillatory stagnation is detected}

\If{ $k+1<mp$}{
    \Return{false}\;
}

\For{$\ell=0,1,\ldots,m-1$}{

Define the $\ell$-th subsequence of $\mathcal{S}^{(k)}$ as
\[
\mathcal{W}_\ell
=
\left(
\alpha_{\,k-mp+\ell p+1},
\ldots,
\alpha_{\,k-mp+(\ell+1)p}
\right).
\]

Compute its range:
\[
R_\ell
=
\max_{\alpha_i\in\mathcal{W}_\ell}\alpha_i
-
\min_{\alpha_i\in\mathcal{W}_\ell}\alpha_i.
\]

Compute its midpoint:
\[
M_\ell
=
\frac{1}{2}
\left(
\max_{\alpha_i\in\mathcal{W}_\ell}\alpha_i
+
\min_{\alpha_i\in\mathcal{W}_\ell}\alpha_i
\right).
\]
}

Compute the variations in the ranges and midpoints:
\[
\Delta_R
=
\max_{0\leq \ell\leq m-1} R_\ell
-
\min_{0\leq \ell\leq m-1} R_\ell,
\]
and
\[
\Delta_M
=
\max_{0\leq \ell\leq m-1} M_\ell
-
\min_{0\leq \ell\leq m-1} M_\ell.
\]

\uIf{$\Delta_R \leq \varepsilon_R$
    \textbf{and}
    $\Delta_M \leq \varepsilon_M$}{
    \Return{true}\;
}
\Else{
    \Return{false}\;
}

\end{minipage}%
}
\end{algorithm}

For each subsequence \(\mathcal{W}_\ell\), the range \(R_\ell\) measures the vertical spread of the convergence values and therefore quantifies the oscillation amplitude, while the midpoint \(M_\ell\) identifies the level around which the oscillation occurs. Small values of both \(\Delta_R\) and \(\Delta_M\) indicate that several consecutive subsequences exhibit similar amplitudes and remain centered around approximately the same convergence level. This behavior indicates that the algorithm is not making sufficient progress toward consensus.

When Algorithm \ref{stagnation_detection} detects oscillatory stagnation, i.e., when the variations in both $R_\ell$ and $M_\ell$ across the considered subsequences fall below their respective tolerances $\varepsilon_R$ and $\varepsilon_M$, the penalty parameter is increased by one half of its initial value:
\begin{equation}
\rho^{(k+1)}
=
\rho^{(k)}
+
\frac{\rho^{(0)}}{2}.
\label{PHA_penalty_update_stagnation}
\end{equation}
Otherwise, the current penalty value is retained. Increasing \(\rho\)
strengthens the penalty imposed on deviations from the consensus and encourages the scenario-specific first-stage solutions to move toward a common decision.

Algorithm \ref{APHA_TS} presents the enhanced version of the classical PHA, which integrates the rounding heuristic and the adaptive penalty update mechanism based on oscillatory stagnation detection Algorithm \ref{stagnation_detection}.

\begin{algorithm}[h!]

\SetKwInOut{KwIn}{Input}
\SetKwInOut{KwOut}{Output}

\scalebox{0.8}{%
\begin{minipage}{1.25\linewidth}

\KwIn{Penalty parameter $\rho^0=|\Omega|$, maximum number of iterations $K$, 
tolerance $\epsilon>0$, subsequence size $p$.}
\KwOut{A near-optimal first-stage solution $\bar{x}^{(k+1)}$.}

\tcc{--- Initialization step ---}
\ForEach{$\omega \in \Omega$ \tcp*[f]{in parallel}}{
    Solve:
    \begin{equation}
    x_{\omega}^{*(0)}
    =
    \arg\max_{x\in \mathcal{X}_{\omega}}
    f(x;\omega)
    \label{InitialStep}
    \end{equation}
}

$\bar{x}^{(0)}
\leftarrow
\sum_{\omega \in \Omega}
p_{\omega}\,x_{\omega}^{*(0)}$\;

$\lambda_{\omega}^{(0)}
\leftarrow
\rho^0
\big(
x_{\omega}^{*(0)}-\bar{x}^{(0)}
\big),
\quad
\forall \omega \in \Omega$\;

$\alpha^{(0)}
\leftarrow
\sum_{\omega \in \Omega}
p_{\omega}
\left\|
x_{\omega}^{*(0)}
-
\bar{x}^{(0)}
\right\|_2$\;

\tcc{--- Main EPHA loop ---}
\For{$k=0,\ldots,K-1$}{

    \ForEach{$\omega \in \Omega$ \tcp*[f]{in parallel}}{

        Solve the augmented subproblem:
        \begin{equation}
        x_{\omega}^{*(k+1)}
        =
        \arg\max_{x\in \mathcal{X}_{\omega}}
        \Big(
        f(x;\omega)
        -
        (\lambda_{\omega}^{(k)})^\top x
        -
        \tfrac{\rho^k}{2}
        \sum_{i\in N}
        \big(1-2\bar{x}_i^{(k)}\big)x_i
        \Big).
        \label{AugmModelPHA}
        \end{equation}
    }

    $\bar{\mu}^{(k+1)}
    \leftarrow
    \sum_{\omega \in \Omega}
    p_{\omega}\,x_{\omega}^{*(k+1)}$\;

    \If{$(k+1)\bmod p = 0
        \;\text{and}\;
        \texttt{Stagnation}
        (k,\mathcal{S}^{(k)})=\texttt{true}$}{
        $\rho^k
        \leftarrow
        \rho^k+\dfrac{\rho^0}{2}$\;
    }

    $\bar{x}^{(k+1)}
    \leftarrow
    \texttt{Rounding}
    (\bar{\mu}^{(k+1)},
    \overline{\kappa},
    \underline{\kappa})$\;

    $\lambda_{\omega}^{(k+1)}
    \leftarrow
    \lambda_{\omega}^{(k)}
    +
    \rho^k
    \big(
    x_{\omega}^{*(k+1)}
    -
    \bar{x}^{(k+1)}
    \big),
    \quad
    \forall \omega \in \Omega$\;

    $\alpha^{(k+1)}
    \leftarrow
    \sum_{\omega \in \Omega}
    p_{\omega}
    \left\|
    x_{\omega}^{*(k+1)}
    -
    \bar{x}^{(k+1)}
    \right\|_2$\;

    \If{$\alpha^{(k+1)} \le \epsilon$}{
        \textbf{break}\;
    }
}

\Return $\bar{x}^{(k+1)}$\;

\end{minipage}%
}

\caption{Enhanced PHA (EPHA) for the TSSQKP.}
\label{APHA_TS}

\end{algorithm}

In Algorithm \ref{APHA_TS}, the condition $(k+1)\bmod p = 0$ ensures that the stagnation detection procedure is applied only at the end of each subsequence of size $p$, while the procedure itself verifies that at least $m$ complete subsequences are available for analysis. Let \(x_{EPHA}\) be the solution obtained from Algorithm \ref{APHA_TS}, and let \(x^{*}\) denote the optimal solution of the original problem. We denote the EPHA lower bound by:
\begin{equation}
	\mathrm{LB}_{\mathrm{EPHA}} = \mathcal{F}(x_{EPHA};\Omega) \leq \mathcal{F}(x^{*};\Omega).
	\label{EPHALB}
\end{equation}
\section{Scenario reduction methods}
\label{Solution_Meth_SRM}
In this section, we propose methods for computing lower and upper bounds for the considered problem by solving reduced versions of the original scenario set. These methods are based on the scenario-clustering strategy developed by \cite{hewitt2022decision}. The proposed methods rely on grouping similar scenarios into clusters and solving the associated clustered subproblems. Here, similarity is measured in the decision space rather than in parameter space. This choice ensures that scenarios are grouped according to their effect on the optimal decisions, since scenarios with similar parameter values may lead to different solutions, while scenarios with different parameters may produce similar decisions.

Let us consider the set of scenarios \(\Omega\) that represents the empirical probability distribution, with associated probabilities \(p_{\omega}\), where \(\omega \in \Omega\).
Suppose that scenario \(\omega \in \Omega\) occurs with certainty. The corresponding optimal solution is given by
\begin{equation}
x^{*}_{\omega}
=
\underset{x \in \mathcal{X}_{\omega}}{\mathrm{argmax}}\; f(x;\omega).
\end{equation}

In practice, however, the realized scenario is unknown. To address this uncertainty, \cite{hewitt2022decision} proposed an opportunity-cost function that measures the loss incurred when the decision associated with scenario \(\omega\) is applied while scenario \(\omega'\) actually occurs:
\begin{equation}
\delta(\omega|\omega')
=
f(x^{*}_{\omega};\omega)
-
f(x^{*}_{\omega'};\omega)
\geq 0,
\label{Distance}
\end{equation}
for \(\omega \neq \omega'\). The inequality in \eqref{Distance} follows from the optimality of \(x^{*}_{\omega}\) for $f(\cdot,\omega)$.

Using the opportunity-cost function $\delta(.,.)$, a symmetric distance measure between scenarios is defined as
\begin{equation}
d(\omega,\omega')
=
\delta(\omega|\omega')
+
\delta(\omega'|\omega),
\qquad
\omega,\omega' \in \Omega.
\end{equation}

The values \(d(\cdot,\cdot)\) are then used to construct a symmetric distance matrix with zero diagonal entries. This matrix serves as the basis for partitioning the scenario set \(\Omega\) into \(m\) clusters,
\[
\mathcal{C}_1 \cup \mathcal{C}_2 \cup \cdots \cup \mathcal{C}_m
=
\Omega,
\]
where each cluster contains scenarios exhibiting similar decision behavior which means that the scenarios generate identical or nearly identical optimal decisions, or that the decision obtained under one scenario remains near-optimal when evaluated under the other scenario. This procedure requires solving $|\Omega|\times|\Omega|$ scenario-based binary subproblems, which can be significantly time consuming. Therefore, we consider a continuous relaxation of these subproblems, which preserves \eqref{Distance} while reducing the computational time required to construct the distance matrix.

Let \(f_{\mathrm{Relax}}(x;\omega)\) denote the relaxed TSSQKP associated with scenario \(\omega \in \Omega\), in which the binary decision variables are relaxed to continuous variables in \([0,1]\). The function \(f_{\mathrm{Relax}}(\cdot,\cdot)\) is used in Algorithm \ref{relax_opp_matrix_math} (see Appendix \ref{Dcomp}) to construct the opportunity-cost matrix.

The \(k\)-medoids clustering algorithm \citep{schubert2019faster} takes as input the distance matrix \(\mathcal{D}\) obtained from Algorithm \ref{relax_opp_matrix_math} (see Appendix \ref{Dcomp}), together with the required number of clusters, and partitions the scenario set \(\Omega\) into $k$ clusters of scenarios.

Unlike $k$-means, $k$-medoids selects representative observations from the original dataset as cluster centers, called medoids \citep{keutchayan2023problem}. By exploiting the symmetry of the distance matrix, only the upper triangular elements need to be computed, and all the subproblems can be solved in parallel.

\subsection{Medoid lower bound}
Let \(\mathcal{C}_1,\ldots,\mathcal{C}_m \subseteq \Omega\) denote the scenario clusters constructed from the opportunity-cost distance matrix \(\mathcal{D}\) by using $k$-medoids clustering. For each cluster \(\mathcal{C}_k\), a representative scenario \(\sigma_k \in \mathcal{C}_k\), referred to as the medoid, is selected for \(k=1,\ldots,m\) and $\Omega_{\sigma}=\{\sigma_k:k=1,\ldots,m\}$ the set of the representative scenarios. The probability associated with cluster \(\mathcal{C}_k\) is defined as
\[
\pi_k=\sum_{\omega\in\mathcal{C}_k}p_{\omega},
\]

where $p_{\omega}=\frac{1}{|\Omega|}$. Using the medoid scenarios and their associated probabilities, a reduced TSSQKP with \(m\ll|\Omega|\) scenarios is constructed. The corresponding solution is obtained by solving
\begin{equation}
\tilde{x}
=
\operatorname*{arg\,max}_{x \in \mathcal{X}_{\Omega_{\sigma}}}
\left\{
\sum_{k=1}^{m}
\pi_k f(x;\sigma_k)
\right\}.
\end{equation}

The solution \(\tilde{x}\) is then evaluated over the original scenario set \(\Omega\). The corresponding medoid lower bound proposed in \cite{hewitt2022decision} is therefore defined as
\begin{equation}
\mathrm{LB}_{\mathrm{M}}
=
\sum_{\omega\in\Omega}
p_{\omega}f(\tilde{x};\omega)
\leq
\mathcal{F}(x^*;\Omega),
\label{MB}
\end{equation}
where \(x^*\) denotes the optimal solution of the original TSSQKP.

The reduced problem contains only \(m\ll|\Omega|\) representative scenarios and can therefore be solved more efficiently than the original problem. In particular, the EPHA (Algorithm \ref{APHA_TS}) may be applied to compute a near-optimal solution for the reduced TSSQKP.

\subsection{Cluster bounds}

Let \(\mathcal{C}_1,\ldots,\mathcal{C}_m \subseteq \Omega\) denote the $m$ clusters obtained from the opportunity-cost distance matrix. For each cluster \(\mathcal{C}_k\), an optimal cluster decision \(x^{*}_{\mathcal{C}_k}\) is computed by solving the following reduced TSSQKP:
\begin{equation}
x^{*}_{\mathcal{C}_k}
=
\underset{x\in \mathcal{X}_{\mathcal{C}_k}}{\mathrm{argmax}}
\;
\sum_{\omega\in \mathcal{C}_k}
p_{\omega} f(x;\omega),
\qquad
k=1,\ldots,m.
\label{SolutionCk}
\end{equation}

Equation \eqref{SolutionCk} defines \(m\) reduced TSSQKPs. Since the clusters may contain different numbers of scenarios, some subproblems can become computationally challenging when \(|\mathcal{C}_k|\) is large. The clustering procedure does not explicitly control the distribution of scenarios among clusters.

To compute the solutions \(x^{*}_{\mathcal{C}_k}\), we employ both a commercial optimization solver and the EPHA described in Section \ref{APHA_TS}, since some cluster subproblems may remain difficult to solve due to their scenario sizes. The resulting cluster solutions are then used to derive lower and upper bounds for the original TSSQKP. 

As mentioned in \cite{hewitt2022decision}, the cluster solutions 
\(x^{*}_{\mathcal{C}_k}\) give rise to the quantities
\begin{align}
\mathrm{LB}_{\mathrm{C}}
&=
\max_{k=1,\ldots,m}
\sum_{\omega \in \Omega}
p_{\omega} f(x^{*}_{\mathcal{C}_k};\omega),
\label{CLB}
\\
\mathrm{UB}_{\mathrm{C}}
&=
\sum_{k=1}^{m}
\sum_{\omega\in \mathcal{C}_k}
p_{\omega} f(x^{*}_{\mathcal{C}_k};\omega),
\label{CUB}
\end{align}
which satisfy
\begin{equation}
\mathrm{LB}_{\mathrm{C}}
\leq
\mathcal{F}(x^{*};\Omega)
\leq
\mathrm{UB}_{\mathrm{C}}.
\label{CBs}
\end{equation}

The cluster lower bound ($\mathrm{LB}_{\mathrm{C}}$) is obtained by evaluating each cluster solution \(x^{*}_{\mathcal{C}_k}\) over the complete scenario set \(\Omega\) and selecting the best resulting feasible solution. The cluster upper bound ($\mathrm{UB}_{\mathrm{C}}$) is obtained by allowing each cluster to use its own optimal decision, which relaxes the nonanticipativity requirement of the original problem and therefore provides an upper bound for the maximization problem.

In the remainder of this work, we refer to the bounds defined in \eqref{MB}, \eqref{CLB}, and \eqref{CUB} as the medoid lower bound ($\mathrm{LB}_{\mathrm{M}}$), the cluster lower bound ($\mathrm{LB}_{\mathrm{C}}$), and the cluster upper bound ($\mathrm{UB}_{\mathrm{C}}$), respectively. When the optimal solution cannot be obtained exactly, the cluster upper bound ($\mathrm{UB}_{\mathrm{C}}$) is used to estimate the optimality gap with respect to the lower bounds $\mathrm{LB}_{\mathrm{EPHA}}$, $\mathrm{LB}_{\mathrm{M}}$, and $\mathrm{LB}_{\mathrm{C}}$, since, from \eqref{CBs}, the true optimal value lies between the lower and upper bounds. When EPHA is used to compute the cluster solutions \(x^{*}_{\mathcal{C}_k}\), the resulting $\mathrm{LB}_{\mathrm{C}}$ remains a valid lower bound provided that the selected solution is feasible for the complete scenario set \(\Omega\), whereas the corresponding $\mathrm{UB}_{\mathrm{C}}$ is only an estimate of the upper bound because EPHA does not guarantee optimal solutions to the cluster subproblems. We refer to this estimated upper bound as the estimated cluster upper bound ($\mathrm{EUB}_{\mathrm{C}}$).

\section{Numerical experiments}
\label{Num_Ex}

In this section, we present the numerical experiments conducted to evaluate the performance of the proposed methods. All methods were implemented in the C programming language using the \texttt{gcc 12.2.0} compiler. The MIP and LP models were solved using CPLEX version 22.1.1 \citep{cplex2022}, while Python version 3.12.3 with the \texttt{kmedoids} package \citep{schubert2022fastkmedoids} was used for the clustering procedure. All computational experiments were performed on a single cluster node of the Digital Research Alliance of Canada supercomputer equipped with a 2.4 GHz processor, providing 64 CPU cores. The scenario subproblems were solved in parallel using 8 CPU cores. The implementation required at most \(100\) GB of RAM.

\subsection{Data description}

We randomly generated instances of the TSSQKP for different profit matrix densities
\[
\Delta \in \{25\%, 50\%, 75\%, 100\%\}.
\]

The density $\Delta$ of an instance refers to the proportion of nonzero quadratic profit terms between pairs of items. The instances were generated according to the scheme proposed in \cite{gallo1980quadratic}. In total, \(800\) instances were generated, with \(200\) instances for each density level. The experiments consider problem sizes $n \in \{10,15,20,25,30\}$, and numbers of scenarios $
|\Omega| \in \{50,80,100,150\}.
$

For each instance, the first-stage profit matrix \(Q\) is generated as a symmetric matrix. 
Its coefficients are positive integers sampled from a uniform distribution in \([1,100]\). 
The first-stage weights \(w_i\) are integers randomly generated in the interval \([1,100]\). 

For the second stage, \(|\Omega|\) scenarios are considered. 
For each scenario \(\omega\), two symmetric profit matrices \(T^{\omega}\) and \(R^{\omega}\) are generated. 
Their coefficients are integers uniformly distributed in \([1,50]\). The scenarios are sampled independently.

The second-stage weights \(a_i^{\omega}\) are integers randomly generated in \([50,100]\). 
The value \(M\) is computed as the minimum scenario weight sum, i.e.,
\[
M=\min_{\omega=1,\ldots,|\Omega| }\sum_{i=1}^{n} a_i^{\omega}.
\]

Finally, the capacity \(\beta\) is sampled uniformly between \(50\) and 
\[
\min\left\{\sum_{i=1}^{n} w_i,\, M\right\},
\]
with
\[
\sum_{i=1}^{n} w_i > 50
\quad \text{and} \quad
M > 50.
\]

The value 
\[
\min\left\{\sum_{i=1}^{n} w_i,\, M\right\}
\]
is used because the capacity is assumed to be identical in both the first and second stages, thereby avoiding scenarios in which the total item weight is smaller than the capacity. The off-diagonal coefficients are generated according to a density parameter \(\Delta\), which represents the percentage of nonzero values in the matrices \(Q\), \(T^{\omega}\), and \(R^{\omega}\). The full data set, including the test instances and optimal solutions for some of them, is available in \cite{dandije2026tssqkp}.

\subsection{Results and analysis}

In this section, we present the computational results. First, we provide a preliminary analysis of the optimal solutions obtained for the TSSQKP. Specifically, we report the optimality gaps with respect to the optimal solution obtained by CPLEX within a time limit of \(5\) hours, together with the lower bounds $\mathrm{LB}_{\mathrm{EPHA}}$, the medoid lower bound ($\mathrm{LB}_{\mathrm{M}}$), and the cluster lower bound ($\mathrm{LB}_{\mathrm{C}}$).

EPHA was executed using the following parameters: a maximum number of iterations \(K = 100\), a convergence tolerance \(\epsilon = 10^{-5}\), and a penalty parameter \(\rho^0\) initially set equal to the total number of scenarios in the problem. The scenario subproblems arising in \eqref{InitialStep} and \eqref{AugmModelPHA} of the EPHA described in Algorithm \ref{APHA_TS}, as well as the cluster-based subproblems defined in \eqref{SolutionCk}, were solved in parallel using \(8\) CPU cores of a single compute node.

For the clustering-based approach, we report the best values of $\mathrm{LB}_{\mathrm{C}}$ and $\mathrm{LB}_{\mathrm{M}}$ obtained using the optimal number of clusters chosen from the set \(\{8,12,16,20\}\). The clusters were generated using the distance matrix computed by Algorithm \ref{relax_opp_matrix_math}. The optimal number of clusters was determined using the silhouette score criterion, a metric for evaluating the quality of a clustering result \citep{rousseeuw1987silhouettes}.

Second, for all instances, we computed the cluster lower bound ($\mathrm{LB}_{\mathrm{C}}$), the cluster upper bound ($\mathrm{UB}_{\mathrm{C}}$), the medoid lower bound ($\mathrm{LB}_{\mathrm{M}}$), and the EPHA lower bound ($\mathrm{LB}_{\mathrm{EPHA}}$), each within a maximum computational time of \(5\) hours. The computational time required to obtain the cluster bounds, $\mathrm{LB}_{\mathrm{C}}$ and $\mathrm{UB}_{\mathrm{C}}$, was then used as the time limit for CPLEX, ensuring a fair comparison between the proposed framework and the solver. The optimality gaps achieved by CPLEX within this time limit were compared with those obtained by the proposed methods. All subproblems arising in the proposed methods were solved in parallel using \(8\) threads. The last set of results presents the optimality gaps obtained by applying EPHA to solve the subproblems associated with each cluster and the corresponding set of medoid scenarios.

\subsubsection{Preliminary results}
\label{PR}

In this section, we compare the performance of CPLEX, EPHA, the cluster lower bound ($\mathrm{LB}_{\mathrm{C}}$), and the medoid lower bound ($\mathrm{LB}_{\mathrm{M}}$) on the instances solved to optimality by CPLEX within a time limit of \(5\) hours. The experiments were conducted for instance densities of \(25\%\), \(50\%\), \(75\%\), and \(100\%\), and for different numbers of scenarios.

The optimality gap associated with each lower bound, namely $\mathrm{LB}_{\mathrm{EPHA}}$, $\mathrm{LB}_{\mathrm{C}}$, and $\mathrm{LB}_{\mathrm{M}}$, is computed as follows:
\begin{equation}
 \text{Gap}(\%) =
 \left(
 \frac{\mathrm{LB}_{\mathrm{CPLEX}} - \mathrm{LB}_{\mathrm{Method}}}
 {\mathrm{LB}_{\mathrm{CPLEX}}}
 \right) \times 100\%.
\end{equation}

Since the methods that provide $\mathrm{LB}_{\mathrm{EPHA}}$, $\mathrm{LB}_{\mathrm{C}}$, and $\mathrm{LB}_{\mathrm{M}}$ yield near-optimal solutions, the above gap evaluates the quality of the obtained lower bounds relative to the optimal solution computed by CPLEX.

Using the CPLEX solver, we were able to solve \(29\) out of \(200\) instances with density \(\Delta = 25\%\), \(13\) out of \(200\) instances with \(\Delta = 50\%\), \(4\) out of \(200\) instances with \(\Delta = 75\%\), and \(6\) out of \(200\) instances with \(\Delta = 100\%\) within the \(5\)-hour time limit. The computational results comparing CPLEX, EPHA, $\mathrm{LB}_{\mathrm{C}}$, and $\mathrm{LB}_{\mathrm{M}}$ are reported in Table \ref{Opt_Gap}.

\begin{table}[htbp]
\centering
\setlength{\tabcolsep}{8pt}
\renewcommand{\arraystretch}{0.7}

\begin{adjustbox}{width=\textwidth}
\begin{tabular}{r rrr r rr rr rr r}
\toprule

\multirow{2}{*}{Density} &
\multirow{2}{*}{${n}$} &
\multirow{2}{*}{${|\Omega|}$} &
\multirow{2}{*}{Instance} &
\multirow{2}{*}{${C^*}$} &
\multicolumn{2}{c}{$\mathrm{LB}_{\mathrm{EPHA}}$ \eqref{EPHALB}} &
\multicolumn{2}{c}{$\mathrm{LB}_{\mathrm{C}}$ \eqref{CLB}} &
\multicolumn{2}{c}{$\mathrm{LB}_{\mathrm{M}}$ \eqref{MB}} &
CPLEX \\

\cmidrule(lr){6-7}
\cmidrule(lr){8-9}
\cmidrule(lr){10-11}

& & & & &
Gap (\%) & Time (s) &
Gap (\%) & Time (s) &
Gap (\%) & Time (s) &
Time (s) \\

\midrule
25 & 10 & 50 & 1 & 16 & 0.00 & 4.57 & 0.00 & 23.27 & 0.00 & 21.79 & 9.19 \\
25 & 10 & 50 & 2 & 8 & 0.00 & 10.47 & 0.00 & 26.86 & 0.00 & 20.81 & 1197.21 \\
25 & 10 & 50 & 3 & 12 & 0.00 & 55.07 & 0.00 & 26.66 & 0.00 & 26.63 & 737.03 \\
25 & 10 & 50 & 4 & 8 & 0.00 & 13.94 & 0.00 & 21.95 & 0.00 & 18.01 & 18.77 \\
25 & 10 & 50 & 5 & 8 & 0.00 & 151.60 & 0.00 & 27.14 & 0.00 & 18.87 & 183.69 \\
25 & 10 & 50 & 6 & 12 & 0.00 & 15.99 & 0.00 & 31.77 & 0.00 & 23.12 & 262.83 \\
25 & 10 & 50 & 7 & 12 & 0.00 & 9.99 & 0.00 & 23.27 & 0.00 & 16.93 & 316.26 \\
25 & 10 & 50 & 8 & 20 & 0.00 & 312.42 & 0.00 & 20.86 & 0.00 & 23.07 & 173.91 \\
25 & 10 & 50 & 9 & 8 & 0.00 & 9.59 & 0.00 & 22.21 & 0.00 & 17.76 & 154.98 \\
25 & 10 & 50 & 10 & 12 & 0.00 & 6.42 & 0.00 & 23.83 & 0.00 & 20.27 & 1824.91 \\
25 & 10 & 80 & 2 & 8 & 0.00 & 12.75 & 0.00 & 52.92 & 0.00 & 48.11 & 91.41 \\
25 & 10 & 80 & 4 & 8 & 0.00 & 51.46 & 0.00 & 215.58 & 0.00 & 47.57 & 4041.61 \\
25 & 10 & 80 & 5 & 20 & 0.00 & 3.36 & 0.00 & 48.58 & 0.00 & 48.35 & 3.93 \\
25 & 10 & 80 & 6 & 12 & 0.00 & 3.57 & 0.00 & 43.76 & 0.00 & 41.19 & 40.98 \\
25 & 10 & 80 & 9 & 16 & 0.00 & 6.52 & 0.00 & 46.98 & 0.00 & 45.19 & 6.25 \\
25 & 10 & 100 & 5 & 8 & 0.00 & 11.01 & 0.00 & 126.94 & 0.09 & 101.27 & 5614.27 \\
25 & 10 & 100 & 6 & 12 & 0.00 & 10.54 & 0.00 & 110.40 & 0.00 & 103.33 & 1281.21 \\
25 & 10 & 100 & 8 & 8 & 0.00 & 13.67 & 0.00 & 149.10 & 0.00 & 100.15 & 11563.67 \\
25 & 10 & 100 & 9 & 12 & 0.00 & 32.95 & 0.00 & 106.48 & 0.01 & 96.89 & 2104.36 \\
25 & 10 & 150 & 1 & 8 & 0.00 & 7.60 & 0.00 & 187.28 & 0.00 & 169.46 & 3675.79 \\
25 & 10 & 150 & 5 & 12 & 0.00 & 18.94 & 0.00 & 206.81 & 0.04 & 200.01 & 3893.36 \\
25 & 10 & 150 & 7 & 20 & 0.00 & 7.67 & 0.00 & 195.96 & 0.00 & 194.19 & 8.38 \\
25 & 10 & 150 & 9 & 16 & 0.00 & 63.88 & 0.00 & 310.43 & 0.00 & 190.40 & 7566.97 \\
25 & 15 & 50 & 2 & 16 & 0.00 & 25.57 & 0.00 & 38.81 & 0.00 & 34.72 & 321.43 \\
25 & 15 & 50 & 4 & 12 & 0.00 & 100.27 & 0.00 & 83.18 & 0.00 & 56.37 & 14002.29 \\
25 & 15 & 50 & 7 & 8 & 0.00 & 84.54 & 0.00 & 61.51 & 0.00 & 46.18 & 2273.66 \\
25 & 15 & 50 & 10 & 8 & 0.00 & 43.10 & 0.00 & 455.95 & 0.00 & 29.64 & 15080.09 \\
25 & 15 & 80 & 9 & 8 & 0.00 & 42.38 & 0.00 & 248.24 & 0.00 & 69.54 & 1670.16 \\
25 & 15 & 100 & 5 & 12 & 0.00 & 16.41 & 0.00 & 125.48 & 0.00 & 109.08 & 493.57 \\
\midrule
\multicolumn{5}{c}{Average} & 0.00 & 39.53 & 0.00 & 105.59 & 0.00 & 66.86 & 2710.76 \\
\midrule
50 & 10 & 50 & 1 & 8 & 0.00 & 56.56 & 0.00 & 40.42 & 0.02 & 18.11 & 3531.16 \\
50 & 10 & 50 & 2 & 8 & 0.00 & 43.63 & 0.00 & 31.43 & 0.00 & 17.48 & 289.44 \\
50 & 10 & 50 & 3 & 8 & 0.00 & 77.36 & 0.00 & 38.46 & 0.02 & 17.58 & 10031.46 \\
50 & 10 & 50 & 5 & 8 & 0.00 & 41.81 & 0.00 & 39.03 & 0.00 & 21.22 & 1921.09 \\
50 & 10 & 50 & 6 & 8 & 0.00 & 198.06 & 0.00 & 37.37 & 0.04 & 20.16 & 731.07 \\
50 & 10 & 50 & 9 & 8 & 0.00 & 22.66 & 0.00 & 35.08 & 0.00 & 18.04 & 147.80 \\
50 & 10 & 50 & 10 & 8 & 0.00 & 20.05 & 0.00 & 25.83 & 0.00 & 15.47 & 84.61 \\
50 & 10 & 80 & 5 & 8 & 0.00 & 11.19 & 0.00 & 50.35 & 0.00 & 41.22 & 162.33 \\
50 & 10 & 80 & 7 & 8 & 0.00 & 29.19 & 0.00 & 202.02 & 0.00 & 41.80 & 12745.06 \\
50 & 10 & 100 & 8 & 8 & 0.00 & 112.14 & 0.00 & 130.81 & 0.00 & 90.26 & 10142.42 \\
50 & 10 & 150 & 2 & 8 & 0.00 & 29.36 & 0.00 & 265.72 & 0.00 & 184.95 & 10888.09 \\
50 & 10 & 150 & 6 & 12 & 0.00 & 39.54 & 0.00 & 216.64 & 0.00 & 187.09 & 3743.86 \\
50 & 15 & 50 & 5 & 8 & 0.00 & 13.74 & 0.00 & 31.96 & 0.00 & 24.58 & 203.53 \\
\midrule
\multicolumn{5}{c}{Average} & 0.00 & 53.48 & 0.00 & 88.08 & 0.01 & 53.69 & 4201.69 \\
\midrule
75 & 10 & 50 & 1 & 20 & 0.00 & 51.71 & 0.00 & 39.40 & 0.00 & 35.16 & 3163.64 \\
75 & 10 & 50 & 7 & 8 & 0.00 & 16.79 & 0.00 & 25.47 & 0.05 & 17.81 & 352.03 \\
75 & 10 & 50 & 10 & 8 & 0.00 & 69.44 & 0.00 & 40.25 & 0.04 & 21.28 & 2145.44 \\
75 & 10 & 100 & 7 & 16 & 0.00 & 140.14 & 0.00 & 135.18 & 0.00 & 87.25 & 11927.17 \\
\midrule
\multicolumn{5}{c}{Average} & 0.00 & 69.52 & 0.00 & 60.07 & 0.02 & 40.38 & 4397.07 \\
\midrule
100 & 10 & 50 & 1 & 8 & 0.00 & 33.30 & 0.00 & 69.21 & 0.00 & 26.86 & 2964.15 \\
100 & 10 & 50 & 3 & 8 & 0.00 & 73.55 & 0.00 & 92.88 & 0.04 & 21.24 & 11216.34 \\
100 & 10 & 50 & 4 & 8 & 0.00 & 52.20 & 0.00 & 133.79 & 0.00 & 20.73 & 1288.82 \\
100 & 10 & 50 & 7 & 8 & 0.00 & 47.66 & 0.00 & 46.81 & 0.00 & 20.78 & 160.84 \\
100 & 10 & 80 & 9 & 8 & 0.00 & 17.50 & 0.00 & 69.76 & 0.00 & 57.13 & 194.91 \\
100 & 10 & 150 & 5 & 16 & 0.00 & 22.50 & 0.00 & 215.10 & 0.11 & 189.98 & 5154.66 \\
\midrule
\multicolumn{5}{c}{Average} & 0.00 & 41.12 & 0.00 & 104.59 & 0.02 & 56.12 & 3496.62 \\
\bottomrule
\end{tabular}
\end{adjustbox}

\caption{Optimality gaps and computational times for EPHA, $\mathrm{LB}_{\mathrm{C}}$, $\mathrm{LB}_{\mathrm{M}}$, and CPLEX.}
\label{Opt_Gap}
\end{table}

These results demonstrate that the proposed methods provide near-optimal solutions with extremely small gaps across all density levels. In particular, $\mathrm{LB}_{\mathrm{EPHA}}$ and $\mathrm{LB}_{\mathrm{C}}$ consistently achieved an average gap of \(0.00\%\), while $\mathrm{LB}_{\mathrm{M}}$ produced very small average gaps ranging between \(0.00\%\) and \(0.02\%\). These results indicate that the proposed approaches are highly effective in approximating the optimal solutions obtained by CPLEX. Moreover, the solution quality remains stable as the density and the number of scenarios increase, highlighting the robustness of the methods.

In terms of computational efficiency, EPHA was generally the fastest approach among the proposed methods. The average computational time of EPHA varied between \(39.53\) and \(69.52\) seconds across all density levels, whereas computing $\mathrm{LB}_{\mathrm{C}}$ and $\mathrm{LB}_{\mathrm{M}}$ required additional computational time due to the clustering procedures and the use of CPLEX to solve the associated subproblems. Nevertheless, both $\mathrm{LB}_{\mathrm{C}}$ and $\mathrm{LB}_{\mathrm{M}}$ remained significantly faster than solving the complete problem directly with CPLEX in most instances. By contrast, CPLEX exhibited larger computational times, with average computational times ranging from \(2710.76\) to \(4397.07\) seconds, and several instances requiring more than \(10{,}000\) seconds. These observations demonstrate that the proposed methods provide an excellent trade-off between solution quality and computational effort, particularly for large scenario sets and dense quadratic profit matrices.

\subsubsection{Evaluation of the proposed bounding framework}
\label{Performance}

This second set of experiments was designed to evaluate the effectiveness of the proposed bounding framework and to compare its performance with that of a commercial optimization solver under equivalent computational conditions.

The quality of the bounds produced by the proposed framework was assessed by computing the optimality gaps obtained from the cluster upper bound ($\mathrm{UB}_{\mathrm{C}}$) and from each of the three lower bounds: the cluster lower bound ($\mathrm{LB}_{\mathrm{C}}$), the medoid lower bound ($\mathrm{LB}_{\mathrm{M}}$), and the EPHA lower bound ($\mathrm{LB}_{\mathrm{EPHA}}$). The evaluation was performed over the full set of generated instances using the optimal number of clusters identified in Section \ref{PR}.

To ensure a fair comparison, each instance was also solved using CPLEX with a time limit equal to the computational time required by the proposed framework to obtain the cluster bounds. The computational times associated with the cluster bounds and the medoid lower bound include both the construction of the opportunity-cost distance matrix and the clustering procedure. The optimality gap reported by CPLEX within this time limit was then used as a reference for assessing the effectiveness of the proposed bounding scheme.

The tightness of each lower bound is measured using the following gap metrics, defined as the relative deviation between the upper bound and the corresponding lower bound:
\begin{align}
 \text{Cluster Gap}(\%) 
 &= 
 \left(
 \frac{\mathrm{UB}_{\mathrm{C}} - \mathrm{LB}_{\mathrm{C}}}
 {\mathrm{UB}_{\mathrm{C}}}
 \right) \times 100\%,
 \label{cluster_gap} \\[6pt]
 \text{Medoid Gap}(\%) 
 &= 
 \left(
 \frac{\mathrm{UB}_{\mathrm{C}} - \mathrm{LB}_{\mathrm{M}}}
 {\mathrm{UB}_{\mathrm{C}}}
 \right) \times 100\%,
 \label{medoid_gap} \\[6pt]
 \text{EPHA Gap}(\%) 
 &= 
 \left(
 \frac{\mathrm{UB}_{\mathrm{C}} - \mathrm{LB}_{\mathrm{EPHA}}}
 {\mathrm{UB}_{\mathrm{C}}}
 \right) \times 100\%.
 \label{apha_gap}
\end{align}

A smaller gap indicates a tighter lower bound and a stronger overall bounding scheme, since all methods are evaluated using the same upper bound. The resulting gaps are then compared with the optimality gap reported by CPLEX.

For these experiments, we were able to solve the cluster subproblems and obtain cluster bounds within the \(5\)-hour time limit for \(124\) out of \(200\) instances with density \(\Delta = 25\%\), \(85\) out of \(200\) instances with \(\Delta = 50\%\), \(57\) out of \(200\) instances with \(\Delta = 75\%\), and \(47\) out of \(200\) instances with \(\Delta = 100\%\). The corresponding results are presented in Table \ref{Bounding_Frame_Work}, which provides a comparative analysis of the optimality gaps obtained by CPLEX and by the proposed bounding approaches for the four density levels \(\Delta \in \{25\%,50\%,75\%,100\%\}\). Overall, the results show that the proposed methods produce smaller and more stable optimality gaps than CPLEX under the same computational time limit.

\begin{table}[htbp]
\centering
\scriptsize
\setlength{\tabcolsep}{6pt}
\renewcommand{\arraystretch}{0.7}
\begin{adjustbox}{width=\textwidth}
\begin{tabular}{rrrrrrrrrrr}
\toprule
\multirow{2}{*}{Density}
& \multirow{2}{*}{$n$}
& \multirow{2}{*}{$|\Omega|$}
& \multirow{2}{*}{\shortstack{Instances Solved\\(\ref{SolutionCk})}}
& \multicolumn{1}{c}{CPLEX}
& \multicolumn{2}{c}{Cluster}
& \multicolumn{2}{c}{Medoid}
& \multicolumn{2}{c}{EPHA} \\
\cmidrule(lr){5-5}
\cmidrule(lr){6-7}
\cmidrule(lr){8-9}
\cmidrule(lr){10-11}
& & & 
& \multicolumn{1}{c}{Gap (\%)}
& \multicolumn{1}{c}{Gap (\%)}
& \multicolumn{1}{c}{Time (s)}
& \multicolumn{1}{c}{Gap (\%)}
& \multicolumn{1}{c}{Time (s)}
& \multicolumn{1}{c}{Gap (\%)}
& \multicolumn{1}{c}{Time (s)} \\
\midrule
\multirow{16}{*}{25} & \multirow{4}{*}{10} & 50 & 10/10 & 5.65 & 4.09 & 24.78 & 4.18 & 20.73 & 4.09 & 59.01 \\
 &  & 80 & 10/10 & 5.37 & 3.99 & 85.77 & 4.73 & 50.84 & 4.00 & 26.28 \\
 &  & 100 & 10/10 & 5.19 & 2.68 & 131.53 & 3.67 & 101.56 & 2.68 & 29.10 \\
 &  & 150 & 10/10 & 7.16 & 2.30 & 261.53 & 2.74 & 201.39 & 2.36 & 58.23 \\
 & \multirow{4}{*}{15} & 50 & 10/10 & 16.30 & 3.50 & 157.95 & 3.77 & 95.99 & 3.50 & 80.98 \\
 &  & 80 & 10/10 & 11.17 & 2.04 & 582.83 & 2.97 & 85.04 & 2.04 & 131.98 \\
 &  & 100 & 10/10 & 7.95 & 2.27 & 1103.14 & 3.29 & 204.91 & 2.28 & 274.11 \\
 &  & 150 & 6/10 & 9.05 & 2.33 & 3771.01 & 3.04 & 401.02 & 2.33 & 234.90 \\
 & \multirow{4}{*}{20} & 50 & 10/10 & 14.05 & 2.19 & 738.65 & 2.56 & 263.73 & 2.19 & 420.14 \\
 &  & 80 & 8/10 & 16.32 & 2.24 & 1730.77 & 2.28 & 416.23 & 2.11 & 537.89 \\
 &  & 100 & 6/10 & 5.90 & 0.56 & 3187.54 & 0.65 & 250.03 & 0.56 & 176.95 \\
 &  & 150 & 6/10 & 5.87 & 0.53 & 6959.74 & 1.00 & 536.08 & 0.53 & 371.86 \\
 & \multirow{4}{*}{25} & 50 & 9/10 & 15.72 & 1.32 & 3450.77 & 1.44 & 877.31 & 1.33 & 5625.64 \\
 &  & 80 & 5/10 & 8.12 & 1.26 & 5724.03 & 1.46 & 750.87 & 1.26 & 1105.24 \\
 &  & 100 & 3/10 & 10.83 & 1.58 & 6265.35 & 1.58 & 2086.06 & 1.58 & 382.34 \\
 &  & 150 & 1/10 & 5.43 & 0.76 & 7332.83 & 0.76 & 1170.32 & 0.84 & 793.43 \\
\midrule
\multirow{15}{*}{50} & \multirow{4}{*}{10} & 50 & 10/10 & 19.20 & 2.49 & 42.92 & 3.24 & 18.99 & 2.49 & 57.66 \\
 &  & 80 & 10/10 & 12.40 & 2.15 & 294.22 & 2.15 & 53.87 & 2.15 & 81.10 \\
 &  & 100 & 10/10 & 13.10 & 1.45 & 395.88 & 1.83 & 105.00 & 1.45 & 158.96 \\
 &  & 150 & 8/10 & 5.88 & 1.51 & 3500.72 & 2.67 & 184.77 & 3.71 & 237.80 \\
 & \multirow{4}{*}{15} & 50 & 10/10 & 34.40 & 2.28 & 831.65 & 3.49 & 110.99 & 2.28 & 335.11 \\
 &  & 80 & 8/10 & 16.88 & 0.85 & 2714.91 & 1.48 & 230.30 & 0.85 & 541.36 \\
 &  & 100 & 6/10 & 9.00 & 0.70 & 5924.52 & 1.43 & 179.82 & 0.72 & 941.41 \\
 &  & 150 & 4/10 & 19.00 & 2.10 & 5181.41 & 2.64 & 1527.91 & 2.10 & 1303.89 \\
 & \multirow{4}{*}{20} & 50 & 9/10 & 25.67 & 1.41 & 1952.74 & 1.82 & 850.16 & 1.41 & 792.87 \\
 &  & 80 & 3/10 & 19.67 & 0.95 & 6880.53 & 0.95 & 798.53 & 0.95 & 2919.10 \\
 &  & 100 & 1/10 & 5.00 & 0.25 & 4008.35 & 0.25 & 15358.89 & 0.25 & 1374.34 \\
 &  & 150 & 1/10 & 37.00 & 3.12 & 11940.56 & 3.12 & 1621.89 & 3.15 & 2462.37 \\
 & \multirow{2}{*}{25} & 50 & 2/10 & 17.00 & 0.75 & 5925.14 & 0.88 & 4310.87 & 0.75 & 3451.79 \\
 &  & 80 & 1/10 & 6.00 & 1.94 & 3365.65 & 5.38 & 158.74 & 1.94 & 241.57 \\
 & \multirow{1}{*}{30} & 50 & 2/10 & 7.50 & 0.87 & 5510.84 & 0.87 & 7051.14 & 0.87 & 6927.93 \\
\midrule
\multirow{10}{*}{75} & \multirow{4}{*}{10} & 50 & 10/10 & 23.50 & 2.42 & 62.43 & 3.53 & 57.76 & 2.42 & 67.08 \\
 &  & 80 & 9/10 & 29.67 & 1.94 & 193.71 & 2.42 & 78.96 & 1.94 & 119.55 \\
 &  & 100 & 9/10 & 40.67 & 2.26 & 429.86 & 3.04 & 103.96 & 2.26 & 211.26 \\
 &  & 150 & 6/10 & 10.83 & 1.54 & 1427.75 & 1.54 & 234.82 & 1.54 & 96.67 \\
 & \multirow{4}{*}{15} & 50 & 7/10 & 29.43 & 2.86 & 1187.24 & 2.86 & 2363.78 & 2.86 & 887.11 \\
 &  & 80 & 7/10 & 28.00 & 1.10 & 2357.93 & 1.30 & 564.52 & 1.10 & 589.35 \\
 &  & 100 & 3/10 & 40.00 & 2.40 & 10789.54 & 3.11 & 575.68 & 2.46 & 2259.60 \\
 &  & 150 & 2/10 & 25.00 & 1.56 & 3658.14 & 1.56 & 552.24 & 1.56 & 473.76 \\
 & \multirow{1}{*}{20} & 50 & 2/10 & 29.50 & 0.44 & 1490.46 & 0.87 & 4978.30 & 0.44 & 2525.34 \\
 & \multirow{1}{*}{25} & 50 & 2/10 & 7.00 & 1.07 & 12745.94 & 2.93 & 449.51 & 1.07 & 3673.93 \\
\midrule
\multirow{10}{*}{100} & \multirow{4}{*}{10} & 50 & 10/10 & 25.80 & 2.39 & 171.73 & 3.03 & 71.54 & 2.39 & 55.54 \\
 &  & 80 & 8/10 & 33.50 & 1.73 & 292.58 & 1.81 & 134.25 & 1.75 & 680.28 \\
 &  & 100 & 9/10 & 30.89 & 3.35 & 3283.60 & 4.40 & 95.09 & 3.35 & 744.08 \\
 &  & 150 & 8/10 & 19.62 & 2.54 & 667.44 & 4.34 & 302.60 & 2.55 & 848.66 \\
 & \multirow{2}{*}{15} & 50 & 5/10 & 39.80 & 2.08 & 2826.00 & 2.12 & 3075.43 & 2.08 & 457.17 \\
 &  & 80 & 3/10 & 30.33 & 1.39 & 4485.89 & 4.29 & 740.09 & 1.39 & 497.75 \\
 & \multirow{3}{*}{20} & 50 & 1/10 & 2.00 & 0.33 & 13460.92 & 0.33 & 471.01 & 0.33 & 2057.82 \\
 &  & 80 & 1/10 & 1.00 & 0.14 & 17105.68 & 0.52 & 457.04 & 0.14 & 473.80 \\
 &  & 100 & 1/10 & 99.00 & 0.70 & 6286.24 & 0.70 & 3446.76 & 0.70 & 757.37 \\
 & \multirow{1}{*}{30} & 50 & 1/10 & 22.00 & 1.38 & 1540.60 & 2.77 & 227.26 & 1.38 & 1393.09 \\
\midrule
\multicolumn{4}{r}{\text{Average}} & \text{18.93} & \text{1.77} & \text{3616.51} & \text{2.35} & \text{1157.93} & \text{1.81} & \text{1000.09} \\
\bottomrule
\end{tabular}
\end{adjustbox}
\caption{First approach bounding gaps and computational times for EPHA, $\mathrm{UB}_{\mathrm{C}}$, $\mathrm{LB}_{\mathrm{C}}$, $\mathrm{LB}_{\mathrm{M}}$, and CPLEX.}
\label{Bounding_Frame_Work}
\end{table}

For all density levels, CPLEX produces larger gaps, with several instances reaching gaps above \(40\%\). In contrast, the gaps associated with $\mathrm{LB}_{\mathrm{C}}$, $\mathrm{LB}_{\mathrm{M}}$, and $\mathrm{LB}_{\mathrm{EPHA}}$ remain concentrated around \(2\%-3\%\), indicating that the proposed bounds remain tight across the tested density levels.

As the density increases, the CPLEX gaps deteriorate considerably. For \(\Delta = 100\%\), the gap approaches \(99\%\) for \(n = 20\) and \(\lvert \Omega \rvert = 50\). These results suggest that the exact optimization approach becomes increasingly difficult to apply efficiently as the quadratic structure becomes denser. By contrast, the gaps associated with $\mathrm{LB}_{\mathrm{M}}$ and $\mathrm{LB}_{\mathrm{EPHA}}$ remain slightly larger than those associated with $\mathrm{LB}_{\mathrm{C}}$ but continue to provide competitive lower bounds with limited dispersion.

Regarding computational time, the computation of $\mathrm{LB}_{\mathrm{M}}$ requires the shortest execution times across all density levels, while EPHA also remains computationally efficient, with moderate variability. The computation of the cluster bounds requires greater computational effort due to the solution of multiple cluster subproblems, each of which is solved using CPLEX.

The computational experiments therefore highlight a balance between solution quality and computational effort. The cluster bounds and $\mathrm{LB}_{\mathrm{EPHA}}$ provide the tightest gaps, whereas $\mathrm{LB}_{\mathrm{M}}$ offers a competitive lower bound with significantly shorter execution times, mainly because it considers only a small number of representative scenarios.

Finally, we observe that, even with scenario clustering, it was not possible to solve all \(200\) instances for every density level within the \(5\)-hour time limit when using CPLEX to solve the cluster subproblems. This limitation is mainly due to the fact that some clusters may still contain a large number of scenarios, leading to difficult sub-TSSQKPs. The next set of experiments therefore investigates the use of EPHA to solve the cluster subproblems for the computation of the cluster bounds and the medoid lower bound ($\mathrm{LB}_{\mathrm{M}}$). This choice is motivated by the strong computational performance and solution quality previously observed for EPHA compared to the CPLEX solver; see Table \ref{Opt_Gap}.

\subsubsection{Cluster-based and medoid solutions with progressive hedging algorithm}

In this section, we evaluate the performance of EPHA for solving the cluster subproblems generated by the proposed scenario-clustering framework. Since solving all cluster subproblems to optimality with CPLEX remains computationally challenging for large and dense instances, EPHA is considered as a heuristic alternative for computing the cluster lower bound ($\mathrm{LB}_{\mathrm{C}}$), the estimated cluster upper bound ($\mathrm{EUB}_{\mathrm{C}}$), and the medoid lower bound ($\mathrm{LB}_{\mathrm{M}}$); see Figure \ref{structure_tssqkp_1} in Appendix \ref{Illustration} for an illustration.

More specifically, EPHA is used to solve the reduced cluster subproblems associated with $\mathrm{LB}_{\mathrm{C}}$, $\mathrm{EUB}_{\mathrm{C}}$, and $\mathrm{LB}_{\mathrm{M}}$. Since EPHA does not guarantee optimal solutions to the cluster subproblems, the resulting upper bound is only an estimate of the true cluster upper bound. Consequently, the optimality gaps computed using $\mathrm{EUB}_{\mathrm{C}}$ are referred to as estimated optimality gaps ($\mathrm{EGap}$).

The objective of this set of experiments is to assess whether the proposed heuristic framework can produce high-quality solutions for a larger number of instances within the imposed computational time limit. The experiments also examine the resulting estimated optimality gaps and computational times when EPHA is integrated into the clustering-based framework. The results therefore provide insight into the balance between computational efficiency and bound quality when the cluster subproblems are solved heuristically rather than exactly. The corresponding results are reported in Table \ref{Bounding_Frame_Work_PHABASED}.

\begin{table}[htbp]
\centering
\scriptsize
\setlength{\tabcolsep}{5pt}
\renewcommand{\arraystretch}{0.7}
\begin{adjustbox}{width=\textwidth}
\begin{tabular}{rrrrrrrrrrr}
\toprule
\multirow{2}{*}{Density}
& \multirow{2}{*}{$n$}
& \multirow{2}{*}{$|\Omega|$}
& \multirow{2}{*}{\shortstack{Instances Solved\\(\ref{SolutionCk})}}
& \multicolumn{1}{c}{CPLEX}
& \multicolumn{2}{c}{Cluster}
& \multicolumn{2}{c}{Medoid}
& \multicolumn{2}{c}{EPHA} \\
\cmidrule(lr){5-5}
\cmidrule(lr){6-7}
\cmidrule(lr){8-9}
\cmidrule(lr){10-11}
& & &
& \multicolumn{1}{c}{Gap (\%)}
& \multicolumn{1}{c}{$\mathrm{EGap}$ (\%)}
& \multicolumn{1}{c}{$\mathrm{EGap}$ (s)}
& \multicolumn{1}{c}{$\mathrm{EGap}$ (\%)}
& \multicolumn{1}{c}{Time (s)}
& \multicolumn{1}{c}{$\mathrm{EGap}$ (\%)}
& \multicolumn{1}{c}{Time (s)} \\
\midrule
\multirow{20}{*}{25} & \multirow{4}{*}{10} & 50 & 10/10 & 4.10 & 3.97 & 65.53 & 4.34 & 25.55 & 3.97 & 59.01 \\
 &  & 80 & 9/10 & 5.78 & 4.41 & 114.84 & 5.24 & 52.92 & 4.42 & 28.80 \\
 &  & 100 & 10/10 & 4.60 & 2.67 & 155.91 & 5.30 & 104.65 & 2.67 & 29.10 \\
 &  & 150 & 8/10 & 6.75 & 2.84 & 331.72 & 3.39 & 202.44 & 2.92 & 70.02 \\
 & \multirow{4}{*}{15} & 50 & 10/10 & 4.10 & 3.29 & 428.02 & 3.64 & 53.22 & 3.27 & 80.98 \\
 &  & 80 & 9/10 & 12.78 & 1.99 & 477.23 & 3.01 & 73.47 & 1.99 & 127.29 \\
 &  & 100 & 10/10 & 7.30 & 2.06 & 787.89 & 3.09 & 124.48 & 2.07 & 274.11 \\
 &  & 150 & 10/10 & 5.30 & 1.74 & 8673.34 & 2.63 & 300.38 & 1.74 & 255.05 \\
 & \multirow{4}{*}{20} & 50 & 10/10 & 4.90 & 1.91 & 2706.49 & 2.27 & 189.93 & 1.91 & 420.14 \\
 &  & 80 & 10/10 & 18.30 & 1.86 & 3167.13 & 2.19 & 155.86 & 1.86 & 575.10 \\
 &  & 100 & 10/10 & 18.60 & 0.92 & 2210.17 & 1.18 & 179.11 & 0.92 & 473.98 \\
 &  & 150 & 10/10 & 34.00 & 0.81 & 1888.91 & 2.00 & 441.67 & 0.82 & 433.43 \\
 & \multirow{4}{*}{25} & 50 & 8/10 & 7.88 & 1.52 & 8614.48 & 1.66 & 247.99 & 1.53 & 6757.31 \\
 &  & 80 & 8/10 & 13.38 & 1.68 & 6061.06 & 3.15 & 359.86 & 1.68 & 1381.49 \\
 &  & 100 & 9/10 & 60.78 & 0.72 & 7877.18 & 1.63 & 361.43 & 0.72 & 1745.60 \\
 &  & 150 & 10/10 & 28.30 & 1.34 & 8936.98 & 2.73 & 1010.58 & 1.35 & 3414.73 \\
 & \multirow{4}{*}{30} & 50 & 8/10 & 11.50 & 1.28 & 7883.44 & 1.91 & 647.58 & 1.07 & 3030.91 \\
 &  & 80 & 4/10 & 8.00 & 2.00 & 13635.22 & 4.15 & 683.51 & 2.18 & 4269.94 \\
 &  & 100 & 7/10 & 25.14 & 0.61 & 9325.36 & 1.91 & 584.05 & 0.61 & 4695.52 \\
 &  & 150 & 4/10 & 157.50 & 1.16 & 9794.33 & 1.43 & 913.81 & 1.16 & 4684.52 \\
\midrule
\multirow{20}{*}{50} & \multirow{4}{*}{10} & 50 & 8/10 & 1.38 & 2.27 & 1888.73 & 2.18 & 30.75 & 2.27 & 44.48 \\
 &  & 80 & 7/10 & 5.00 & 2.60 & 1991.60 & 2.68 & 48.17 & 2.60 & 104.91 \\
 &  & 100 & 7/10 & 3.29 & 1.91 & 2421.02 & 2.14 & 115.92 & 1.91 & 174.12 \\
 &  & 150 & 8/10 & 4.12 & 1.50 & 2105.47 & 1.69 & 192.90 & 3.70 & 237.80 \\
 & \multirow{4}{*}{15} & 50 & 10/10 & 16.30 & 2.23 & 1188.34 & 3.43 & 77.82 & 2.23 & 335.11 \\
 &  & 80 & 10/10 & 20.50 & 0.98 & 1179.01 & 1.55 & 91.88 & 0.98 & 515.78 \\
 &  & 100 & 8/10 & 18.75 & 1.52 & 1886.94 & 2.07 & 153.33 & 1.54 & 1049.58 \\
 &  & 150 & 9/10 & 19.22 & 2.11 & 3314.51 & 2.61 & 412.67 & 2.18 & 1245.31 \\
 & \multirow{4}{*}{20} & 50 & 10/10 & 14.90 & 1.21 & 3229.97 & 1.59 & 144.84 & 1.21 & 769.74 \\
 &  & 80 & 7/10 & 21.57 & 0.60 & 4255.19 & 1.26 & 271.02 & 0.60 & 1631.64 \\
 &  & 100 & 8/10 & 38.62 & 1.53 & 8265.99 & 2.59 & 282.57 & 2.70 & 4980.96 \\
 &  & 150 & 7/10 & 37.29 & 2.10 & 8181.58 & 2.78 & 668.91 & 2.11 & 3681.49 \\
 & \multirow{4}{*}{25} & 50 & 3/10 & 17.00 & 0.21 & 14087.77 & 0.29 & 578.30 & 0.21 & 2873.94 \\
 &  & 80 & 5/10 & 38.60 & 1.44 & 11526.46 & 3.45 & 485.51 & 1.44 & 4674.21 \\
 &  & 100 & 5/10 & 23.60 & 0.65 & 10130.10 & 0.65 & 384.82 & 0.65 & 4083.96 \\
 &  & 150 & 4/10 & 79.75 & 1.21 & 9683.01 & 1.52 & 814.18 & 1.21 & 6623.05 \\
 & \multirow{4}{*}{30} & 50 & 2/10 & 2.00 & 0.12 & 10035.10 & 0.12 & 1050.47 & 0.12 & 6927.93 \\
 &  & 80 & 1/10 & 46.00 & 1.28 & 8349.16 & 1.28 & 315.67 & 1.28 & 3198.74 \\
 &  & 100 & 1/10 & 612.00 & 0.16 & 1350.23 & 0.16 & 328.78 & 0.16 & 1649.17 \\
 &  & 150 & 2/10 & 70.50 & 0.69 & 8386.21 & 1.88 & 962.50 & 0.69 & 1653.72 \\
\midrule
\multirow{18}{*}{75} & \multirow{4}{*}{10} & 50 & 9/10 & 10.89 & 2.68 & 217.38 & 3.92 & 45.04 & 2.68 & 71.92 \\
 &  & 80 & 7/10 & 13.14 & 2.04 & 518.46 & 3.09 & 78.32 & 2.04 & 142.45 \\
 &  & 100 & 9/10 & 22.00 & 2.13 & 396.77 & 2.31 & 109.98 & 2.13 & 210.97 \\
 &  & 150 & 10/10 & 25.00 & 1.75 & 561.27 & 1.75 & 225.23 & 1.75 & 140.86 \\
 & \multirow{4}{*}{15} & 50 & 6/10 & 19.33 & 1.68 & 1406.77 & 2.79 & 113.99 & 1.68 & 332.81 \\
 &  & 80 & 7/10 & 35.14 & 1.06 & 1390.23 & 1.16 & 94.23 & 1.06 & 673.05 \\
 &  & 100 & 8/10 & 44.75 & 3.21 & 3205.70 & 3.86 & 263.68 & 3.27 & 3354.15 \\
 &  & 150 & 8/10 & 32.25 & 2.29 & 3977.88 & 3.32 & 410.98 & 2.39 & 3251.79 \\
 & \multirow{4}{*}{20} & 50 & 5/10 & 17.60 & 1.60 & 12659.49 & 1.88 & 1374.68 & 1.60 & 5198.33 \\
 &  & 80 & 6/10 & 28.33 & 1.18 & 9136.75 & 1.65 & 465.84 & 1.18 & 3927.07 \\
 &  & 100 & 6/10 & 46.83 & 1.66 & 13104.72 & 2.73 & 741.64 & 1.66 & 4923.98 \\
 &  & 150 & 5/10 & 51.40 & 1.63 & 14691.20 & 2.61 & 926.46 & 1.63 & 8041.48 \\
 & \multirow{4}{*}{25} & 50 & 4/10 & 31.00 & 0.36 & 7848.87 & 1.46 & 498.78 & 0.36 & 5721.54 \\
 &  & 80 & 2/10 & 21.50 & 1.35 & 3514.49 & 1.35 & 298.38 & 1.35 & 1998.14 \\
 &  & 100 & 4/10 & 39.25 & 0.83 & 6159.93 & 0.83 & 377.18 & 0.84 & 1901.44 \\
 &  & 150 & 4/10 & 50.75 & 1.59 & 13020.40 & 1.81 & 1059.47 & 1.59 & 7338.31 \\
 & \multirow{2}{*}{30} & 100 & 2/10 & 46.50 & 0.05 & 8627.22 & 0.09 & 825.88 & 10.89 & 10862.18 \\
 &  & 150 & 1/10 & 46.00 & 0.47 & 8263.07 & 2.02 & 973.39 & 0.47 & 3368.49 \\
\midrule
\multirow{18}{*}{100} & \multirow{4}{*}{10} & 50 & 7/10 & 23.14 & 3.34 & 197.57 & 4.26 & 29.90 & 3.34 & 60.31 \\
 &  & 80 & 8/10 & 10.25 & 1.94 & 398.48 & 2.79 & 71.06 & 1.96 & 684.34 \\
 &  & 100 & 8/10 & 42.00 & 3.46 & 573.28 & 4.17 & 90.72 & 3.46 & 820.85 \\
 &  & 150 & 6/10 & 17.83 & 2.60 & 900.94 & 4.30 & 257.14 & 2.62 & 1089.71 \\
 & \multirow{4}{*}{15} & 50 & 8/10 & 38.75 & 2.29 & 1345.32 & 2.60 & 128.34 & 2.29 & 477.60 \\
 &  & 80 & 8/10 & 35.00 & 1.18 & 2617.79 & 2.72 & 131.79 & 1.18 & 616.81 \\
 &  & 100 & 8/10 & 26.50 & 1.65 & 3175.34 & 1.75 & 232.62 & 1.65 & 1292.97 \\
 &  & 150 & 8/10 & 50.62 & 1.72 & 7097.43 & 2.22 & 395.10 & 1.76 & 3644.89 \\
 & \multirow{4}{*}{20} & 50 & 8/10 & 34.57 & 0.79 & 12072.61 & 0.79 & 2187.30 & 3.33 & 7854.48 \\
 &  & 80 & 10/10 & 44.90 & 1.61 & 7872.70 & 2.08 & 863.81 & 1.61 & 6135.30 \\
 &  & 100 & 6/10 & 57.00 & 2.06 & 14311.96 & 3.11 & 1310.65 & 2.27 & 8670.46 \\
 &  & 150 & 2/10 & 72.50 & 0.94 & 6653.54 & 0.94 & 641.03 & 0.94 & 3183.75 \\
 & \multirow{4}{*}{25} & 50 & 1/10 & 64.00 & 1.12 & 1100.36 & 1.12 & 367.06 & 1.12 & 1737.12 \\
 &  & 80 & 2/10 & 32.50 & 2.49 & 6066.70 & 3.97 & 477.57 & 2.49 & 3357.20 \\
 &  & 100 & 1/10 & 20.00 & 0.20 & 2103.57 & 0.20 & 256.52 & 0.21 & 1005.41 \\
 &  & 150 & 1/10 & 26.00 & 3.31 & 2933.27 & 5.88 & 844.68 & 3.31 & 749.41 \\
 & \multirow{2}{*}{30} & 50 & 1/10 & 5.00 & 1.38 & 3112.24 & 2.77 & 186.87 & 1.38 & 1393.09 \\
 &  & 80 & 1/10 & 14.00 & 0.66 & 4091.03 & 6.30 & 508.39 & 0.66 & 1876.24 \\
\midrule
\multicolumn{4}{r}{\text{Average}} & \text{35.86} & \text{1.65} & \text{5156.79} & \text{2.41} & \text{420.88} & \text{1.88} & \text{2439.10} \\
\bottomrule
\end{tabular}
\end{adjustbox}
\caption{Second approach estimated bounding gaps and computational times for $\mathrm{EUB}_{\mathrm{C}}$, $\mathrm{LB}_{\mathrm{EPHA}}$, $\mathrm{LB}_{\mathrm{C}}$, $\mathrm{LB}_{\mathrm{M}}$, and CPLEX.}
\label{Bounding_Frame_Work_PHABASED}
\end{table}

Compared with the previous CPLEX implementation of the clustering framework, this variant allows a larger number of instances to be processed within the \(5\)-hour time limit. In particular, cluster bounds are obtained for \(174\) out of \(200\) instances with \(\Delta = 25\%\), \(122\) out of \(200\) instances with \(\Delta = 50\%\), \(103\) out of \(200\) instances with \(\Delta = 75\%\), and \(94\) out of \(200\) instances with \(\Delta = 100\%\). This confirms that replacing the exact solution of the cluster subproblems with EPHA improves the practical applicability of the clustering-based framework.

The results show that CPLEX continues to produce significantly larger optimality gaps than the estimated gaps associated with the proposed bounding approaches under the same computational time. On average, the CPLEX gap is \(35.86\%\), whereas the average estimated gaps associated with $\mathrm{LB}_{\mathrm{C}}$, $\mathrm{LB}_{\mathrm{M}}$, and $\mathrm{LB}_{\mathrm{EPHA}}$ are only \(1.65\%\), \(2.41\%\), and \(1.88\%\), respectively. Thus, the proposed approaches remain much tighter than CPLEX, even when EPHA is used heuristically to solve the cluster subproblems. Although EPHA solves the cluster subproblems heuristically, the resulting average estimated gaps can be smaller because they are computed using the estimated cluster upper bound ($\mathrm{EUB}_{\mathrm{C}}$), which is not guaranteed to be a valid upper bound on the optimal objective value.

For all density levels, the estimated gaps associated with $\mathrm{LB}_{\mathrm{C}}$ provide the smallest gaps overall. The estimated gaps associated with $\mathrm{LB}_{\mathrm{EPHA}}$ are also highly competitive and remain close to those associated with $\mathrm{LB}_{\mathrm{C}}$ in most cases. The estimated gaps associated with $\mathrm{LB}_{\mathrm{M}}$ are slightly larger on average, but their performance remains stable and better than that of CPLEX. These results confirm that the proposed bounding framework preserves good solution quality even when exact optimization of the cluster subproblems is replaced by a heuristic procedure.

As the density of the quadratic profit matrix increases, the CPLEX gaps remain large and unstable, with several instances exhibiting gaps above \(40\%\). This behavior is particularly pronounced for some larger instances, such as the case \(\Delta = 50\%\), \(n = 30\), and \(\lvert \Omega \rvert = 100\), where the CPLEX gap reaches \(612\%\). In contrast, the estimated gaps associated with $\mathrm{LB}_{\mathrm{C}}$ and $\mathrm{LB}_{\mathrm{EPHA}}$ remain mostly concentrated around \(1\%-3\%\), indicating that the proposed bounds are robust with respect to the density of the profit matrix. Although a few isolated cases exhibit larger estimated gaps associated with $\mathrm{LB}_{\mathrm{EPHA}}$, such as \(10.89\%\) for \(\Delta = 75\%\), \(n = 30\), and \(\lvert \Omega \rvert = 100\), the overall behavior remains stable across the tested density levels.

Regarding computational time, the computation of $\mathrm{LB}_{\mathrm{M}}$ remains the most efficient in terms of time, with an average time of \(420.88\) seconds. EPHA requires more computational effort, with an average time of \(2439.10\) seconds, but remains faster than the computation of the cluster bounds, whose average time is \(5156.79\) seconds. This difference is expected, since the computation of the cluster bounds requires solving several cluster subproblems, whereas $\mathrm{LB}_{\mathrm{M}}$ relies only on representative scenarios.

Overall, the EPHA framework provides a favorable balance between bound quality and computational efficiency. Compared with the CPLEX solution of the cluster subproblems, it increases the number of instances for which bounds can be obtained within the time limit while maintaining tight estimated optimality gaps. The estimated gaps associated with $\mathrm{LB}_{\mathrm{C}}$ provide the best solution quality, $\mathrm{LB}_{\mathrm{EPHA}}$ offers estimated gaps that are very close to those associated with $\mathrm{LB}_{\mathrm{C}}$ with reduced computational effort, and $\mathrm{LB}_{\mathrm{M}}$ remains the fastest alternative while still producing competitive lower bounds. Even though EPHA is a heuristic method, its first iteration may still exceed the \(5\)-hour time limit for large instances, since the scenario subproblems must be solved exactly with CPLEX. This explains why we are not able to solve all instances within the imposed time limit, especially as the number of items increases; see the augmented scenario-based model in Appendix \ref{AugModel}, which is another challenging combinatorial problem.

\section{Conclusion}
\label{Concl}

In this work, we investigated a two-stage stochastic quadratic knapsack problem under uncertainty in profits and weights. To address the computational challenges induced by the stochastic and quadratic structure of the problem, we proposed an enhanced progressive hedging framework combined with clustering-based scenario reduction techniques for computing lower and upper bounds.

First, an enhanced progressive hedging algorithm (EPHA) was developed by integrating a rounding procedure and a dynamic penalty update strategy based on the detection of oscillation and stagnation patterns.

Second, we introduced a clustering framework based on opportunity-cost distances between scenarios. Using the resulting scenario partitions, we derived the medoid lower bound ($\mathrm{LB}_{\mathrm{M}}$), the cluster lower bound ($\mathrm{LB}_{\mathrm{C}}$), and the cluster upper bound ($\mathrm{UB}_{\mathrm{C}}$). This procedure decomposes the original stochastic problem into smaller and easier-to-solve subproblems.

To test the proposed methods, we conducted computational experiments on \(800\) generated TSSQKP instances with different numbers of items, scenario sizes, and profit matrix densities. When the cluster subproblems were solved exactly with CPLEX, the resulting bounds produced substantially smaller optimality gaps than those obtained by CPLEX under comparable computational conditions. In particular, the gaps associated with $\mathrm{LB}_{\mathrm{C}}$ were generally the smallest across the tested density levels, while $\mathrm{LB}_{\mathrm{M}}$ and $\mathrm{LB}_{\mathrm{EPHA}}$ also provided competitive lower bounds.

The experiments also showed that the performance of CPLEX deteriorates as the density of the quadratic profit matrix increases, whereas the proposed framework remains relatively stable. This suggests that the clustering-based strategy scales more favorably with problem density than directly solving the complete problem with CPLEX.

From a computational perspective, the computation of $\mathrm{LB}_{\mathrm{M}}$ required the shortest computational times, while the cluster bounds provided the best bound quality at the cost of greater computational effort. The experiments also showed that solving all cluster subproblems exactly remains challenging for large and dense instances, since some clusters may still contain many scenarios and lead to difficult combinatorial subproblems.

Motivated by these observations, EPHA was also investigated as a heuristic approach for solving the cluster subproblems. In this case, $\mathrm{LB}_{\mathrm{C}}$ remains a valid lower bound provided that the selected cluster solution is feasible for the complete scenario set, whereas the corresponding cluster upper bound is no longer guaranteed to be valid because EPHA does not necessarily solve the cluster subproblems to optimality. We therefore denote the resulting quantity by the estimated cluster upper bound ($\mathrm{EUB}_{\mathrm{C}}$), and the corresponding gaps are interpreted as estimated optimality gaps. This EPHA-based variant increases the number of instances for which bounds can be obtained within the \(5\)-hour time limit while maintaining small estimated gaps. However, EPHA may still fail to solve large instances efficiently, since the scenario subproblems themselves become increasingly challenging as the number of items grows.

Future research may focus on developing dedicated heuristic procedures for solving the individual scenario subproblems arising within EPHA, particularly for large-scale instances where the quadratic knapsack structure becomes increasingly difficult as the number of items grows. Another promising direction is the design of branch-and-cut or decomposition-based exact algorithms for solving the cluster subproblems more efficiently, together with the investigation of alternative scenario similarity measures.

\section*{Acknowledgements}
We thank the Digital Research Alliance of Canada for providing high-performance computing facilities. 
\section*{Declarations}
\subsection*{Data availability}  
Enquiries about data availability should be directed to the authors
\subsection*{Funding}  
This work was partly supported by the Canadian Natural Sciences and Engineering Research Council under grants 2021-03307 and 2025-03964. This support is greatly appreciated.
\subsection*{Competing Interests}  
The authors declare that they have no competing interests.

%%%%%%%%%%%%%%%%%%%%%%%

%%%%%%%%%%%%%%%%

\section{Appendix}
\subsection{Augmented model}
\label{AugModel}
For a given scenario $\omega\in\Omega$, in Algorihm \ref{APHA_TS} at iteration $k$, we solve the following augmented model:

\begin{center}
\begin{minipage}{0.95\textwidth}
\begingroup

\setlength{\jot}{-3pt}
\setlength{\abovedisplayskip}{2pt}
\setlength{\belowdisplayskip}{2pt}
\setlength{\abovedisplayshortskip}{2pt}
\setlength{\belowdisplayshortskip}{2pt}

\begin{subequations}
\label{APHA_Augmented_Model}

%--------------------------------------------------------
% Objective function
%--------------------------------------------------------
\begin{equation}
\resizebox{0.98\linewidth}{!}{$
\displaystyle
\max_{x,y,u,v,z,t}\quad
\sum_{i\in N} Q_{ii}x_i
+
2\sum_{(i,j)\in\mathcal{P}} Q_{ij}y_{ij}
+
\sum_{i\in N} T^{\omega}_{ii}u^\omega_i
+
2\sum_{(i,j)\in\mathcal{P}} T^\omega_{ij}z^\omega_{ij}
-
\sum_{i\in N} R^\omega_{ii}v^\omega_i
-
2\sum_{(i,j)\in\mathcal{P}} R^\omega_{ij}t^\omega_{ij}
-
(\lambda_{\omega}^{(k)})^\top x
-
\tfrac{1}{2}\rho^k
\sum_{i\in N}\bigl(1-2\bar{x}_i^{(k)}\bigr)x_i
$}
\label{APHA_Augmented_Obj}
\end{equation}

\vspace{-5pt}

%--------------------------------------------------------
% Constraints
%--------------------------------------------------------
\footnotesize

\begin{alignat}{2}
\text{s.t.}\quad
&
\sum_{i\in N} w_i x_i \leq \beta,
&&
\label{APHA_Augmented_FS_Capacity}
\\
&
\sum_{i\in N} a_i^\omega
\left(
u_i^\omega+x_i-v_i^\omega
\right)
\leq \beta,
&\quad&
\label{APHA_Augmented_SS_Capacity}
\\
&
y_{ij}\leq x_i,
&\quad&
\forall (i,j)\in\mathcal{P},
\label{APHA_Augmented_Y1}
\\
&
y_{ij}\leq x_j,
&\quad&
\forall (i,j)\in\mathcal{P},
\label{APHA_Augmented_Y2}
\\
&
y_{ij}\geq x_i+x_j-1,
&\quad&
\forall (i,j)\in\mathcal{P},
\label{APHA_Augmented_Y3}
\\
&
u_i^\omega\leq 1-x_i,
&\quad&
\forall i\in N,
\label{APHA_Augmented_Addition}
\\
&
v_i^\omega\leq x_i,
&\quad&
\forall i\in N,
\label{APHA_Augmented_Removal}
\\
&
z_{ij}^\omega\leq u_i^\omega,
&\quad&
\forall (i,j)\in\mathcal{P},
\label{APHA_Augmented_Z1}
\\
&
z_{ij}^\omega\leq u_j^\omega,
&\quad&
\forall (i,j)\in\mathcal{P},
\label{APHA_Augmented_Z2}
\\
&
z_{ij}^\omega
\geq u_i^\omega+u_j^\omega-1,
&\quad&
\forall (i,j)\in\mathcal{P},
\label{APHA_Augmented_Z3}
\\
&
t_{ij}^\omega\leq v_i^\omega,
&\quad&
\forall (i,j)\in\mathcal{P},
\label{APHA_Augmented_T1}
\\
&
t_{ij}^\omega\leq v_j^\omega,
&\quad&
\forall (i,j)\in\mathcal{P},
\label{APHA_Augmented_T2}
\\
&
t_{ij}^\omega
\geq v_i^\omega+v_j^\omega-1,
&\quad&
\forall (i,j)\in\mathcal{P},
\label{APHA_Augmented_T3}
\\
&
x_i,u_i^\omega,v_i^\omega\in\{0,1\},
&\quad&
\forall i\in N,
\label{APHA_Augmented_Binary1}
\\
&
y_{ij},z_{ij}^\omega,t_{ij}^\omega\in\{0,1\},
&\quad&
\forall (i,j)\in\mathcal{P}.
\label{APHA_Augmented_Binary2}
\end{alignat}

\end{subequations}

\endgroup
\end{minipage}
\end{center}

% \begin{equation}
% \resizebox{0.8\textwidth}{!}{$
% \begin{aligned}
% \max_{x,u,v}\;& \sum_{i\in N} Q_{ii}x_i
% +\sum_{i\in N}\sum_{\substack{j\in N \\ j\neq i}} Q_{ij}y_{ij} \\
% &+
% \sum_{i\in N} T^{\omega}_{ii}u^\omega_i
% +\sum_{i\in N}\sum_{\substack{j\in N \\ j\neq i}} T^\omega_{ij}z^\omega_{ij}
% -\sum_{i\in N} R^\omega_{ii}v^\omega_i
% -\sum_{i\in N}\sum_{\substack{j\in N \\ j\neq i}} R^\omega_{ij}t^\omega_{ij}\\
% &-(\lambda_{\omega}^{(k)})^\top x
% 			-\tfrac{1}{2}\rho^k
% 			\sum_{i\in N}\big(1-2\bar{x}_i^{(k)}\big)x_i
% \\[4pt]
% \text{s.t. } &
% \sum_{i\in N} w_ix_i \le \beta,\\
% &
% \sum_{i\in N} a^\omega_i(u^\omega_i+x_i-v^\omega_i)
% \le \beta,\\
% &
% y_{ij}\le x_i,
% \qquad \forall (i<j)\in N,\\
% &
% y_{ij}\le x_j,
% \qquad \forall (i<j)\in N,\\
% &
% y_{ij}\ge x_i+x_j-1,
% \qquad \forall (i<j)\in N,\\
% &
% u^\omega_i \le 1-x_i,
% \qquad \forall i\in N,\\
% &
% v^\omega_i \le x_i,
% \qquad \forall i\in N,\\
% &
% z^\omega_{ij}\le u^\omega_i,
% \qquad \forall (i<j)\in N,\\
% &
% z^\omega_{ij}\le u^\omega_j,
% \qquad \forall (i<j)\in N,\\
% &
% z^\omega_{ij}\ge u^\omega_i+u^\omega_j-1,
% \qquad \forall (i<j)\in N,\\
% &
% t^\omega_{ij}\le v^\omega_i,
% \qquad \forall (i<j)\in N,\\
% &
% t^\omega_{ij}\le v^\omega_j,
% \qquad \forall (i<j)\in N,\\
% &
% t^\omega_{ij}\ge v^\omega_i+v^\omega_j-1,
% \qquad \forall (i<j)\in N,\\
% &
% x_i,y_{ij},u^\omega_i,v^\omega_i,z^\omega_{ij},t^\omega_{ij}
% \in \{0,1\}.
% \end{aligned}
% $}
% \label{TSSQKP_Linearise}
% \end{equation}

\subsection{Distance matrix computation}
\label{Dcomp}
\begin{algorithm}[H]
	\caption{Opportunity-cost distance matrix construction}
	\label{relax_opp_matrix_math}
	
	\DontPrintSemicolon
	
	\SetKwInOut{Input}{Input}
	\SetKwInOut{Output}{Output}
	
	\Input{Instance of TSSQKP.}
	
	\Output{Opportunity-cost distance matrix \(\mathcal{D}\).}
	
	\BlankLine
	
	\textbf{Step 1: Solve relaxed scenario subproblems independently.}\;
	
	For each scenario \(\omega \in \Omega\), assume that the scenario occurs with certainty and solve:
	\[
	x^{R,*}_{\omega}
	=
	\underset{x}{\mathrm{argmax}}\;
	f_{\mathrm{Relax}}(x;\omega).
	\]
	
	\BlankLine
	
	\textbf{Step 2: Compute pairwise opportunity-cost distances.}\;
	
	For every pair of scenarios \((\omega,\omega')\), with \(\omega,\omega' \in \Omega\) and $\omega\ne \omega'$, compute
	\begin{equation*}
	d(\omega,\omega')
	=
	\delta(\omega|\omega')
	+
	\delta(\omega'|\omega),
	\end{equation*}
	where
	\begin{equation*}
	\delta(\omega|\omega')
	=
	f_{\mathrm{Relax}}(x^{R,*}_{\omega'};\omega')
	-
	f_{\mathrm{Relax}}(x^{R,*}_{\omega};\omega').
	\end{equation*}
	
	Set
	\[
	d(\omega',\omega)=d(\omega,\omega')
	\quad \text{and} \quad
	d(\omega,\omega)=0.
	\]
	
	\BlankLine
	
	\textbf{Step 3: Construct the opportunity-cost matrix.}\;
	
	Form the matrix
	\[
	\mathcal{D}
	=
	\big[
	d(\omega,\omega')
	\big]_{\omega\ne \omega'}.
	\]
	
\end{algorithm}

\subsection{Numerical experiment schemes}
\label{Illustration}

\begin{figure}[htbp]
	\centering
	\includegraphics[width=0.95\linewidth]{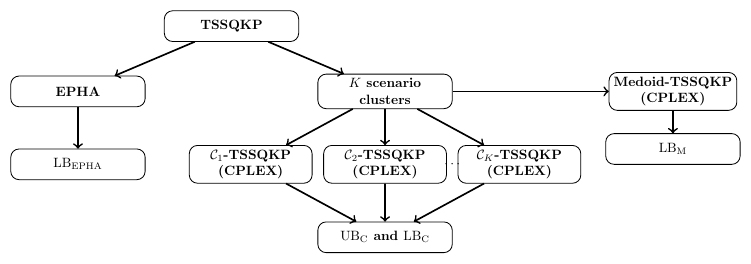}
	\caption{First approach for solving the TSSQKP.}
	\label{structure_tssqkp}
\end{figure}

\begin{figure}[htbp]
	\centering
	\includegraphics[width=0.95\linewidth]{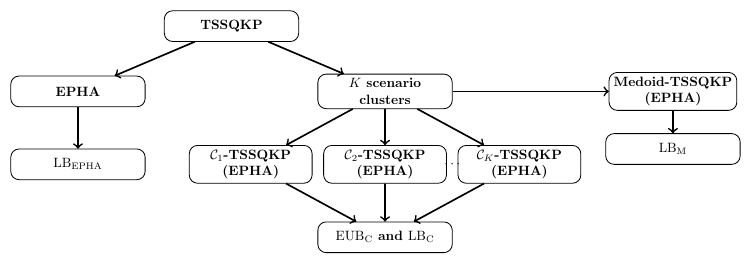}
	\caption{Second approach for solving the TSSQKP.}
	\label{structure_tssqkp_1}
\end{figure}

\end{document}